\documentclass[11pt]{article}
\usepackage{amsmath,amscd,amsfonts,latexsym}
\usepackage{theorem}

\theoremstyle{plain}

\newtheorem{theorem}{Theorem}[section]
\newtheorem{claim}{Claim}[theorem]
\newtheorem{proposition}[theorem]{Proposition}
\newtheorem{lemma}[theorem]{Lemma}
\newtheorem{corollary}[theorem]{Corollary}

\newtheorem{definition}[theorem]{Definition}

\newtheorem{remark}[theorem]{Remark}
\newtheorem{conjecture}{Conjecture}[section]

\theorembodyfont{\upshape}
\newtheorem{proof}{Proof}

\newcommand{\bproof}{\begin{proof}}
\newcommand{\eproof}{\end{proof}}

\newcommand{\qed}{\hspace*{\fill} $\Box$} 

\newcommand{\qedclaim}{\relax}

\newcommand{\cl}[1]{\overline{#1}}

\newcommand{\lca}[1]{{\mathfrak{\lowercase{#1}}}}

\newcommand{\NN}{{\mathbb N}}
\newcommand{\RR}{{\mathbb R}}
\newcommand{\ZZ}{{\mathbb Z}}
\newcommand{\QQ}{{\mathbb Q}}

\def\cU{{\cal U}}
\def\cH{{\cal H}}
 
\newcommand{\inv}{{^{-1}}}

\mathchardef\GG="321D  

\newcommand{\Aut}{\operatorname{Aut}}

\newcommand{\GL}{\operatorname{GL}}
\newcommand{\SL}{{\operatorname{SL}}}
\newcommand{\Ad}{\operatorname{Ad}}

\newcommand{\supp}{\operatorname{supp}}
\newcommand{\Zcl}{\operatorname{Zcl}}
\newcommand{\vol}{\operatorname{vol}}
\newcommand{\rank}{\operatorname{rank}}
\newcommand{\Rad}{\operatorname{Rad}}

\newcommand{\Nsu}{{N^\ast(H_u,u)}}
\newcommand{\capu}{{\cap_{u\in\cU} N(H_u,u)}}
\newcommand{\capsu}{{\cap_{u\in\cU} N^\ast(H_u,u)}}

\title{Invariant Measures and Orbit Closures on Homogeneous Spaces
for Actions of Subgroups Generated by Unipotent Elements}
\author{Nimish A. Shah}

\begin{document}
\maketitle

\begin{abstract}
The theorems of M.~Ratner, describing the finite ergodic invariant
measures and the orbit closures for unipotent flows on homogeneous
spaces of Lie groups, are extended for actions of subgroups generated
by unipotent elements. More precisely: Let $G$ be a Lie group (not
necessarily connected) and $\Gamma$ a closed subgroup of $G$. Let $W$
be a subgroup of $G$ such that $\mbox{Ad}_G(W)$ is contained in the
Zariski closure (in $\mbox{Aut}(\mbox{Lie }G)$) of the subgroup
generated by the unipotent elements of $\mbox{Ad}_G(W)$. Then any
finite $W$-invariant $W$-ergodic measure on $G/\Gamma$ is a
homogeneous measure (i.e.,\ it is supported on a closed orbit of a
subgroup preserving the measure). Moreover, if $G/\Gamma$ has finite
volume (i.e., has a finite $G$-invariant measure), then the closure of
any orbit of $W$ on $G/\Gamma$ is a homogeneous set (i.e.,\ a finite
volume closed orbit of a subgroup containing $W$). Both the above
results hold if $W$ is replaced by any subgroup $\Lambda\subset W$
such that $W/\Lambda$ has finite volume.
\end{abstract}

\section{Introduction}
In \cite{Ratner:measure,Ratner:distribution} Ratner showed the
validity of Raghunathan's conjecture [4] describing orbit closures for
actions of unipotent subgroups on homogeneous spaces of Lie groups, and
its analogous conjecture, due to Dani~\cite{Dani:invent2}, describing
ergodic invariant measures for such actions. Earlier in
\cite{Margulis:varna,Margulis:selberg} Margulis had conjectured that
the conclusions of the orbit closure and the ergodic invariant measure
conjectures should hold also for the actions of subgroups generated by
unipotent elements, as compared to the subgroups themselves being
unipotent. In fact for actions of {\em connected\/} subgroups
 generated by unipotent elements, this conjecture was also verified to
be true in Ratner's above mentioned papers. Using Ratner's theorems
for  actions of unipotent one-parameter subgroups, in this article we
show the validity of the generalized conjecture. This also answers a
question raised by Ratner in \cite[End of Section 4]{Ratner:ICM}.

\bigskip\noindent{\it Notation.} Let $G$ be a Lie group, $\lca{G}$
its Lie algebra, and $\Ad_G:G\to\GL(\lca{G})$ denote the Adjoint
representation of $G$ on $\lca{G}$. An element $u\in G$ is called {\em
$\Ad_G$-unipotent}, if $\Ad_G(u)$ is a unipotent linear
transformation. A subgroup of $G$ consisting of $\Ad_G$-unipotent
elements is called an {\em $\Ad_G$-unipotent subgroup}.

Let $\langle \cU \rangle$ denote the subgroup generated by a subset
$\cU$ in $G$. Let $\Zcl(X)$ denote the Zariski closure of a subset $X$
in $\GL(\lca{G})$. For a subgroup $F$ of $G$, let $F^0$ denote the
connected component of $F$ containing the identity element. 

For a Borel measure $\mu$ on a second countable topological space, we
denote by $\supp(\mu)$ the closed subset which is the complement of the
union of all open sets with zero $\mu$-measure.

\begin{theorem} \label{nimish:thm:measure}
Let $G$ be a Lie group and $\Gamma$ a closed subgroup of $G$. Let $W$
be a subgroup of $G$ and $\cU\subset W$ such that $\cU$ consists of
$\Ad_G$-unipotent elements and $\Ad_G(W)\subset
\Zcl(\Ad_G(\langle\cU\rangle))$. Let $\mu$ be a finite $W$-invariant
$W$-ergodic Borel measure on $G/\Gamma$.  Then there exists a closed
subgroup $H$ of $G$ containing $W$ such that $\mu$ is $H$-invariant and
$\supp(\mu)$ is a closed $H$-orbit. 
\end{theorem}

A Borel measure on a locally compact second countable topological
space is called {\em locally finite}, if it is finite on compact sets.

\begin{theorem} \label{nimish:thm:locally:finite}
Let $G$, $\Gamma$, and $W$ be as in Theorem~\ref{nimish:thm:measure}. Suppose
that $G/\Gamma$ has a finite $G$-invariant measure. Let $\mu$ be a
locally finite $W$-invariant $W$-ergodic measure on $G/\Gamma$. Then
there exists a closed subgroup $H$ of $G$ containing $W$ such that
$\mu$ is $H$-invariant and $\supp(\mu)$ is a closed $H$-orbit.
\end{theorem}

\begin{theorem} \label{nimish:thm:closure}
Let $G$, $\Gamma$ and $W$ be as in Theorem~\ref{nimish:thm:measure}. Suppose
that $G/\Gamma$ has a finite $G$-invariant measure. Then for any $x\in
G/\Gamma$, there exists a closed subgroup $F$ of $G$ containing $W$
such that 
\[
\cl{Wx}=Fx.
\]
Moreover, $F^0x$ has a finite $F^0$-invariant measure (cf.\ Conjectures 
\ref{conj:finite} and \ref{conj:SLn} below). Also the action
of $W$ is ergodic with respect to a locally finite $F$-invariant
measure on $Fx$.
\end{theorem}

We may note that Theorem~\ref{nimish:thm:measure} and
Theorem~\ref{nimish:thm:closure} have already been proved in the above
mentioned papers of Ratner in the following special case: $G$ is
connected, $W$ is of the form $W=\cup_{i=1}^\infty w_i W^0$, where
$w_i$ is $\Ad_G$-unipotent, $i=1,2,\dots$, $W/W^0$ is finitely
generated, and $W^0$ is generated by one-parameter $\Ad_G$-unipotent
subgroups contained in $W^0$. In the case when $G$ is not connected
and $W$ is a nilpotent $\Ad_G$-unipotent subgroup of $G$,
Theorem~\ref{nimish:thm:measure} was proved by Witte
\cite[Theorem~1.2]{Witte:quotients}. In the case when $G$ is connected
and $W$ is a $\Ad_G$-unipotent subgroup, it was shown by
Dani~\cite[Theorem~4.3]{Dani:rk=1} that if $G/\Gamma$ has a finite
invariant measure then any locally finite $W$-invariant $W$-ergodic
measure is finite. In \cite[Remarks~3.12]{Margulis:varna}, Margulis
observed that the same holds for connected $W$. Thus for connected $W$,
Theorem~\ref{nimish:thm:locally:finite} reduces to Theorem~\ref{nimish:thm:measure},
which was proved by Ratner (for connected $W$).

The following result is deduced from that above results using the   
`suspension techniques' (cf.~Witte~\cite[Corollary~5.8]{Witte:quotients}).

\begin{corollary}  \label{nimish:cor:lattice}
Let $G$ and $W$ be as in any one of the theorems stated above. Assume
that $W$ is closed and let $\Lambda$ be a closed subgroup of $W$ such
that $W/\Lambda$ has a finite $W$-invariant measure. Then all the 
theorems stated above are true for $\Lambda$ in place of $W$.

Further, if $G/\Gamma$ admits a finite $G$-invariant measure and $W$
is  connected, then we have the following additional information:
\begin{enumerate}
\item[{\rm 1.}]
Any locally finite $\Lambda$-invariant $\Lambda$-ergodic measure
on $G/\Gamma$ is finite. 
\item[{\rm 2.}] For  $x\in G/\Gamma$, if $\cl{\Lambda x}=Fx$ for a closed 
subgroup $F$ of $G$ then $Fx$ has a finite $F$-invariant measure.
\end{enumerate}
\end{corollary}

From this corollary, we deduce the following.

\begin{corollary} \label{nimish:cor:Lambda.g.Gamma}
Let $G$ be a connected semisimple Lie group without nontrivial compact
factors. Let $\Gamma$ and $\Lambda$ be lattices in $G$ such that at
least one of them is irreducible in $G$; (see \cite[Sect.\
5.20]{Raghunathan:book} for definition). Then either $\Lambda \Gamma$
is dense in $G$ or $\Lambda\cap\Gamma$ is a subgroup of finite index
in $\Gamma$, as well as $\Lambda$.
\end{corollary}

In view of the above results we may ask if the following is true. 

\begin{conjecture} \label{conj:finite} 
Let the notation be as in Theorem~\ref{nimish:thm:measure}. Suppose further
that $G/\Gamma$ has a finite $G$-invariant measure. Then the following
statements hold:
\begin{enumerate}
\item[{\rm 1.}] 
Any locally finite $W$-invariant $W$-ergodic measure on $G/\Gamma$ is
finite.
\item[{\rm 2.}]
For any $x\in G/\Gamma$, if $\cl{Wx}=Fx$ for a closed subgroup $F$ of $G$, 
then $Fx$ has a finite $F$-invariant measure.
\item[{\rm 3.}]
The closure of any $W$-orbit has finitely many connected components.
\end{enumerate}
\end{conjecture}

Note that by the above stated theorems and by Hedlund's
lemma~\ref{lemma:Hedlund}, the three statements in the above
conjecture are equivalent.

\begin{remark} \label{rem:lattice:finite}
If the above conjecture is valid for the diagonal action of $W$ on
$W/\Lambda\times G/\Gamma$, then it holds for the action of $\Lambda$
on $G/\Gamma$, where $W$ and $\Lambda$ are as in
corollary~\ref{nimish:cor:lattice}.
\end{remark}

It seems that the generalized Raghunathan conjecture due to Margulis
already includes Conjecture~\ref{conj:finite}. Using some standard
arguments, as in proof of Theorem~\ref{nimish:thm:closed:finite}, one can
reduce this conjecture to the case of $G$ being a semisimple group
with no nontrivial compact factors and trivial center. Then one can
express $G$ as a product of semisimple subgroups each intersecting
$\Gamma$ in an irreducible lattice. Using the structure of the cusps
in the quotient of the $\RR$-rank one factors, one can take care of
those factors. Thus the conjecture remains to be proved for higher
rank semisimple groups $G$. We use the arithmeticity theorem of
Margulis, and reduce the conjecture to its following typical case.

\begin{conjecture}\!\!\footnote{Recently Alex Eskin and G. A. Margulis informed the author
that they can prove this conjecture.}\label{conj:SLn}%
\hspace*{.5em} Let $G=\SL_n(\RR)$, $\Gamma=\SL_n(\ZZ)$, and $W\subset\SL_n(\QQ)$ a
closed subgroup of $G$ such that $W$ is contained in the Zariski
closure of a subgroup generated by $\Ad_G$-unipotent elements of
$W$. If $W\Gamma$ is discrete, then $W\cap\Gamma$ is of finite index
in $W$.
\end{conjecture}

Finally, a challenging question is to describe invariant measures and
orbit closures for actions of subgroups $H$ whose Zariski closure is
generated by unipotent elements, which are not necessarily contained
in $H$. For example, let $H$ be a Zariski dense subgroup of
$\SL_2(\RR)$ not containing any unipotent elements, and consider the
action of $H$ on $\SL_2(\RR)/\SL_2(\ZZ)$.

\section{Preliminary Results}

In this section we recall some standard  results or their 
modifications about orbit closures on homogeneous spaces, and Zariski
density of certain discrete subgroups as in Borel's density theorem.

\begin{lemma}[Hedlund's Lemma] \label{lemma:Hedlund} 
Let $X$ be a second countable topological space and $W$ be a group of
homeomorphisms of $X$. Let $\mu$ be a $W$-invariant $W$-ergodic Borel
measure on $X$. Then $\cl{Wx}=\supp\mu$ for $\mu$-almost all $x\in X$.
\end{lemma}

\begin{lemma} \label{lemma:finitevol:closed}
Let $G$ be a locally compact second countable group, $\Gamma$ a 
discrete subgroup of $G$, and $\pi:G\to G/\Gamma$ the natural quotient
map. Let $F$ be a Borel measurable subgroup of $G$. Suppose there
exists a locally finite $F$-invariant Borel measure concentrated on
$\pi(F)$. Then $F$ is closed, $\pi(F)$ is closed, and
$\supp(\mu)=\pi(F)$.
\end{lemma}

\proof Since $\mu$ is locally finite, by dominated convergence
theorem (see \cite[Proposition~1.4]{Ratner:solvable}), $\mu$ is
invariant under the closure of $F$ in $G$, say $H$. We have a natural
inclusion $H/H\cap\Gamma\hookrightarrow G/\Gamma$, which is
$H$-equivariant. Since $\mu$ is concentrated on $\pi(F)\subset\pi(H)$,
we can treat $\mu$ as a locally finite $H$-invariant Borel measure on
$H/H\cap\Gamma$. Since $\mu$ is concentrated on an orbit of $F$, we
conclude that a Haar measure on $H$ is strictly positive on $F$. Since
$F=FF^\inv$, we have that $F$ is an open, and hence a closed subgroup
of $H$. Thus $F=H$.

Now since $F$ is closed, the result follows from the proof of
\cite[Theorems~1.12-1.13]{Raghunathan:book}
(cf.~\cite[Proposition~1.4]{Ratner:measure}). Although these
references assume that $\mu$ is finite, the local nature of the
conclusion requires only the assumption that $\mu$ is locally finite.
\qed

\begin{lemma} \label{lemma:centralizer}
Let $G$ be a Lie group and  $\Gamma$ be a discrete subgroup of $G$. Let 
$\pi:G\to G/\Gamma$ be the quotient map. Let $F$ be a subgroup of
$G$ such that $$\Ad_G(F)\subset
\Zcl(\Ad_G(F\cap\Gamma)).$$ Then $\pi(Z_G(F))$ is closed.
\end{lemma}

\proof Take any $\gamma\in\Gamma$, then 
$Z_G(\gamma)\Gamma$ is the inverse image of the discrete set
$\{\lambda\inv \gamma \lambda:\lambda\in \Gamma\}\subset\Gamma$ under
the continuous map $G\ni g\mapsto g\inv \gamma g\in G$. Hence
$\pi(Z_G(\gamma))$ is dense in $G$. Note that if $F_1$ and $F_2$ are
closed subgroups of $G$ such that $\pi(F_1)$ and $\pi(F_2)$ are
closed, then $\pi(F_1\cap F_2)$ is closed. Therefore we conclude that
$\pi(Z_G(F\cap\Gamma))$ is closed. 

By our Zariski closure hypothesis, $Z_G(F\cap\Gamma)^0\subset
Z_G(F)$. Therefore the result follows from an observation that for any
closed subgroup $H$ of $G$, if $\pi(H)$ is closed then $\pi(H_1)$ is
closed for any subgroup $H_1$ of $H$ containing $H^0$.   
\qed

\bigskip\noindent{\bf Definition}
Let $F$ be a connected subgroup of a Lie group $G$ and $\lca{F}$ the
Lie algebra associated to $F$. Let $N_G$ denote the
normalizer of $F$ in $G$. We define
\[
N_G^1(F)=\{g\in N_G(F): \det(\Ad_G(g)|_\lca{F})=1\}.
\]

\begin{remark} \label{rem:nor}
\rm
Note that all $\Ad_G$-unipotent elements of $N_G(F)$ are  contained in
$N_G^1(F)$. Now suppose $W\subset G$ and $\cU\subset W$ such that $\cU$
consists of $\Ad_G$-unipotent elements and $\Ad_G(W)\subset
\Zcl(\Ad_G(\cU))$. Then $$\cU\subset N_G(F)\implies  W\subset
N_G^1(F).$$
\end{remark}
 
\begin{proposition}  \label{nimish:prop:nor}
Let $G$ be a Lie group, $\Gamma$ a discrete subgroup of $G$, and
$\pi:G\to G/\Gamma$ the natural quotient map. Let $U$ be a 
subgroup of $G$ generated by one-parameter $\Ad_G$-unipotent subgroups
of $U$. Suppose $F$ is a closed connected subgroup of $G$ such that
$\pi(F)$ has a finite $F$-invariant measure and 
$\cl{\pi(U)}=\pi(F)$. Then $\pi(N^1(F))$ is closed in $G/\Gamma$.
\end{proposition}

\proof Let $\lca{G}$ be the Lie algebra of $G$ and $\lca{F}$ the Lie
algebra associated to $F$. Let $d=\dim \lca{F}$. Consider the action of
$G$ on $\wedge^d \lca{G}$ via the $d$-exterior power of $\Ad_G$. Let
$p\in \wedge^d \lca{F}\setminus \{0\}$. Then by
\cite[Theorem~3.4]{Dani:Margulis:distribution}, the orbit $\Gamma\cdot
p$ is closed (in the reference it is assumed that $G$ is connected,
but their proof is valid without this assumption).

Observe that the stabilizer of $p$ in $G$ is $N_G^1(F)$. Therefore
$\Gamma N_G^1(F)$ is a closed subset of $G$, and hence the same holds
for $N_G^1(F)\Gamma$.  
\qed



\begin{proposition}[Dani] \label{nimish:prop:closed:discrete}
Let $G$ be a Lie group, $\Gamma$ a closed subgroup of $G$, and
$\pi:G\to G/\Gamma$ the natural quotient map. Let $u\in G$ be an
$\Ad_G$-unipotent element and $\mu$ a finite $u$-invariant measure on
$G/\Gamma$. Then 
\[
\Ad_G(u)\in\Zcl(\Ad_G(g\Gamma g\inv)), \qquad \forall
g\in\pi\inv(\supp(\mu)).
\] 

In particular, if $H$ is a closed subgroup of $G$ containing $u$ such that
$\supp(\mu)\subset\pi(H)$, then 
\[
\Ad_G(u)\in \Zcl(\Ad_G(h(H\cap\Gamma)h\inv), \qquad \forall
h\in H\cap \pi\inv(\supp(\mu)).
\]
\end{proposition}

\proof
This follows from Dani's version of Borel's density
theorem~\cite[Corollary 2.6]{Dani:version:borel:density} (see
\cite[Proof of Corollary~4.3]{Witte:quotients} for details). 
\qed

\section{Extension of a Discrete Unipotent Flow to a Continuous 
Unipotent Flow} \label{sec:extend}

\noindent{\it Notation.} 
Let $G$ be a Lie group and $\Gamma$ be a discrete subgroup of $G$ such
that $G=G^0\Gamma$. Let $\pi:G\to G/\Gamma$ be the natural quotient
map and $x_0=\pi(e)$. Let $u\in G$ be an $\Ad_G$-unipotent element and 
$\lca{G}$ be the Lie algebra of $G$. 

Let $\rho:G_0\to G^0$ be the universal covering homomorphism. Let
$\{\tilde u(t)\}$ be the one-parameter subgroup of $\Aut(G_0)$ such
that 
$$
\{D(\tilde u(t)) \mid_{T_e(G_0)=\lca{G}}\}_{t\in \RR}
=\Zcl(\langle \Ad_G u \rangle)\subset
\Aut(\lca{G}) \mbox{~and~} D(\tilde u(1))\mid_{\lca{G}} =\Ad_G u.
$$ 
Consider the semidirect product ${\bar G}=\RR\cdot G_0$, where $t\in\RR$ acts as
$\tilde u(t)$ on $G_0$; in other words, $t g (-t)= \tilde u(t)(g)$ for
all $g\in G_0$. Note that 
$$
\rho(1g(-1))=u\rho(g)u\inv, \ \forall g\in G_0.
$$ 
Therefore we can extend $\rho:\ZZ\cdot G_0\to \langle u \rangle
G^0$ such that $\rho(1)=u$. Let $\Gamma_1=\rho\inv(\Gamma\cap \langle
u \rangle G^0)$. Since $G=G^0\Gamma$, we have
\begin{equation} \label{eq:identify}
G/\Gamma\cong\langle u\rangle G^0/(\Gamma\cap \langle u \rangle G^0)
\cong \ZZ\cdot G_0/\Gamma_1\subset {\bar G}/\Gamma_1.
\end{equation}
Under this identification, the action of $u$ on $G/\Gamma$ and
the action of $u(1)$ on $\ZZ\cdot G_0/\Gamma_1$ are identical, where 
$u(t)=t\in \bar{G}$ for all $t\in \RR$ and $\{u(t)\}$ is a one-parameter 
$\Ad_{\bar G}$-unipotent subgroup of $\bar G$.

Thus we can treat a discrete unipotent flow as a restriction of a 
continuous unipotent flow. Now we will deduce the algebraic properties
of the invariant measures and orbit closures for the discrete
unipotent flows using the analogous properties of the continuous
unipotent flows.

Let $\mu$ be a finite $u$-invariant $u$-ergodic Borel measure on
$G/\Gamma$. By Hedlund's lemma~\ref{lemma:Hedlund}, there exists $g\in
G^0$ such that $\supp(\mu)=\cl{\langle u \rangle g x_0}$. Let $w=g\inv
u g$ and $\lambda=g\inv\mu$; where by definition, $g\inv\mu(E)=\mu(gE)$
for any Borel subset $E\subset G/\Gamma$. Then $\lambda$ is
$w$-invariant, $w$-ergodic, and $\supp(\lambda)=\cl{\pi(\langle w
\rangle)}$. Let $\tilde g\in \rho\inv(g)$. Put $w(t)= \tilde g\inv
u(t) \tilde g$. We can treat $\lambda$ as a Borel measure on
${\bar G}/\Gamma_1$. Note that the action of $w$ on $G/\Gamma$ and the
action of $w(1)$ on $\ZZ\cdot G_0/\Gamma_1\subset {\bar G}/\Gamma_1$
are isomorphic. Let ${\bar \lambda}$ be the measure on ${\bar
G}/\Gamma_1$ such that for any compactly supported continuous function
$f$ on ${\bar G}/\Gamma_1$, we have
\begin{equation} \label{eq:extend}
\int_{{\bar G}/\Gamma_1} f \, d{\bar \lambda} = \int_0^1\left( 
\int_{{\bar G}/\Gamma_1} f(w(t)x)\, d\lambda(x)\right)\,dt.
\end{equation}
Then ${\bar \lambda}$ is finite, $\{w(t)\}$-invariant, and
$\{w(t)\}$-ergodic. Therefore by  Ratner's measure classification
theorem~\cite[Theorem~1]{Ratner:measure}, there exists a closed
connected subgroup ${\bar H}$ of ${\bar G}$ containing $\{w(t)\}$ such
that ${\bar \lambda}$ is ${\bar H}$-invariant and $\supp({\bar
\lambda})={\bar H}x_0$. Put $H=\rho(\ZZ\cdot G_0\cap {\bar H})$. Then
$H$ is a closed subgroup of $G$ containing $w$ and $\lambda$ is a
finite $H$-invariant measure on $Hx_0$. Therefore by
Lemma~\ref{lemma:finitevol:closed}, $\pi(H)$ is closed,
$\supp(\lambda)=\pi(H)$, and $H\cap\Gamma$ is a lattice in $H$. 

We shall describe orbit closures under the assumption that {\em
$\Gamma$ is a lattice in $G$}. Let $g\in G_0$ and $Y=\cl{\langle
u\rangle gx_0}$.  Let $Z=g\inv Y$ and $w=g\inv u g$. Then
$Z=\cl{\langle w \rangle x_0}$. Let $\tilde g\in \rho\inv(g)$ and
$w(t)=\tilde g\inv u(t) \tilde g$. In view of
Equation~\ref{eq:identify}, $Z=\cl{\langle w(1)\rangle x_0}$. Put
$\tilde Z=w([0,1])Z$. Then $\tilde Z=\cl{\{w(t)\}x_0}\subset \bar
G/\Gamma_1$. By Ratner's description of orbit closures of continuous
unipotent flows~\cite{Ratner:distribution} the following holds: There
exists a closed connected subgroup ${\bar H}$ of ${\bar G}$  containing  
$\{w(t)\}$ such that ${\bar Z}={\bar H}x_0$ and ${\bar Z}$ has a
finite ${\bar H}$-invariant Borel measure, say ${\bar \lambda}$. Also
the trajectory $\{w(t)x_0:t\geq 0\}$ is uniformly  distributed with
respect to $\bar\lambda$. Put $H=\rho(\ZZ\cdot G_0 \cap {\bar
H})$. Then $H$ is a closed subgroup of $G$ containing $w$ such that
$Z=\pi(H)$, and $Z$ has a finite $H$-invariant Borel measure, say
$\lambda$. Also the trajectory $\{w^nx_0:n>0\}$ is uniformly
distributed with respect to $\lambda$.

%
%

\bigskip\noindent{\bf Definition} Let the notation be as in the
beginning of this section. Let $\cH_u$ be the collection of subgroups
$H$ of $G$ such that $H=\langle w \rangle H^0$, $H\cap\Gamma$ is a
lattice in $H$, and $\cl{\langle w\rangle x_0}=Hx_0$, where $w:=g\inv
u g\in H$ for some $g\in G^0$. Let $\lambda_H$ denote a unique
$H$-invariant Borel probability measure on $Hx_0$, for all
$H\in\cH_u$. Note that $\pi(H)$ has finitely many connected components. 

\bigskip
In view of the above discussion and the definitions, we have the
following results: 

\begin{theorem}[Ratner] \label{nimish:thm:ratner:measure}
Let the notation be as in the beginning of this section. Let $\mu$ be a
$u$-invariant $u$-ergodic Borel probability measure on
$G/\Gamma$. Then there exists $g\in G^0$ and $H\in\cH_u$ such that
$ug\in gH$ and $\mu=g\lambda_H$.
\end{theorem}

\begin{theorem}[Ratner] \label{nimish:thm:ratner:closure} 
Let the notation be as in the beginning of this section. Further assume
that $\Gamma$ is a lattice in $G$. Let $g\in G^0$. Then there exists
$H\in\cH_u$ such that $ug\in gH$ and $\cl{\langle u\rangle
\pi(g)}=g\pi(H)$. Moreover, the trajectory $\{u^n\pi(g):n>0\}$ is
uniformly distributed with respect to $g\lambda_H$.
\end{theorem}

\begin{proposition}[Ratner]  \label{nimish:prop:countable}
The collection $\cH_u$ is countable.
\end{proposition} 

\proof 
Let $\cH$ is the collection of all closed connected subgroups ${\bar
H}$ of ${\bar G}$ such that ${\bar H}\cap \Gamma_1$ is a lattice in
${\bar H}$ and for a one-parameter $\Ad_{{\bar G}}$-unipotent
subgroup, say $\{w(t)\}\subset {\bar H}$, we have
$\cl{\{w(t)\}x_0}={\bar H}x_0$. Then by
Proposition~\ref{nimish:prop:closed:discrete} and the countability theorem of
Ratner~\cite[Theorem 1]{Ratner:measure} (see
\cite[Proposition 2.1]{Dani:Margulis:distribution} for another proof), $\cH$
is countable. {}From the above discussion $\cH_u=\{\rho(\ZZ\cdot
G_0\cap {\bar H}):{\bar H}\in\cH\}$. Hence $\cH_u$ is countable.  \qed

\section{Singular Subsets of $G$ Associated to the \\ $u$-action}
\label{sec:singular:subsets}

\noindent{\it Notation.} 
Let $G$ be a Lie group and $\Gamma$ be a discrete subgroup of $G$ such
that $G=G^0\Gamma$. Let $\pi:G\to G/\Gamma$ be the natural quotient
map and $x_0=\pi(e)$. Let $u\in G$ be an $\Ad_G$-unipotent element.

\bigskip\noindent{\bf Definition.}
For $H\in\cH_u$, we say that $F<H$ (or $H>F$) if and only if
$F\in\cH_u$, $F\subset H$, and $\pi(F)\neq\pi(H)$. If $F<H$, then 
either $\dim F< \dim H$, or the number of connected components of 
$\pi(F)$ is less than the  number of connected components of $\pi(H)$. 
Therefore any decreasing sequence $H>F_1>F_2>\cdots$ is finite.

For $H\in\cH_u$, define
\begin{eqnarray*}
N(H,u)&=&\{g\in G^0: ug\in gH\},\\
S(H,u)&=&\bigcup_{F<H}N(F,u) \ \ \ \mbox{~and~}\\ 
N^\ast(H,u)&=&N(H,u)\setminus S(H,u).
\end{eqnarray*}
Note that for any $\gamma\in\Gamma$, $$N(H,u)\gamma=N(\gamma\inv
H\gamma,u) \text{ and } S(H,u)\gamma=S(\gamma\inv H \gamma, u).$$

\begin{proposition} \label{nimish:prop:g:in:N*(H,u)}
Let $H\in\cH_u$ and $g\in N(H,u)$. Then 
\[
g\in N^\ast(H,u)\iff \cl{\langle u \rangle \pi(g)}=g\pi(H).
\] 
\end{proposition}

\proof Replacing $u$ by $g\inv u g$, without loss of generality
we may assume that $g=e$. Since $H\cap\Gamma$ is
a lattice in $H$, by Theorem~\ref{nimish:thm:ratner:closure}, there exists
$F\subset H$ such that $F\in\cH_u$ and $\cl{\pi(\langle u
\rangle)}=\pi(F)$. Now by definition, $$e\in N^\ast(H,u) \iff
\pi(F)=\pi(H).$$ Clearly, $$\pi(F)=\pi(H) \iff \cl{\pi(\langle u
\rangle)}=\pi(H).$$ This completes the proof of the proposition.  \qed

\begin{proposition} \label{nimish:prop:measure:tube}
Let $\lambda$ be a $u$-invariant $u$-ergodic Borel probability  measure on
$G/\Gamma$. Then there exist $H\in\cH_u$ and $g\in N^\ast(H,u)$ such that
$\lambda=g\lambda_H$, where $\lambda_H$
denotes a unique $H$-invariant Borel probability measure on
$\pi(H)$.
\end{proposition}

\proof By Theorem~\ref{nimish:thm:ratner:measure}, there exist $H\in\cH_u$
and $g_1\in N(H,u)$ such that $\lambda=g_1\lambda_H$ and
$\supp\mu=g_1\pi(H)$. By Hedlund's lemma, there exists $h\in H$ such
that $\cl{\langle u\rangle \pi(g_1h)}=\supp(\mu)$. Put $g=g_1h$. Then
$g\in N(H,u)$, $\lambda=g\lambda_H$, and $\cl{\langle u \rangle
\pi(g)}=g\pi(H)$. Now the proposition follows from
Proposition~\ref{nimish:prop:g:in:N*(H,u)}.  \qed

\begin{proposition} \label{nimish:prop:N*(H,u)}
Suppose $g\in N^\ast(H,u)$ and $\gamma\in\Gamma$ such that $g\gamma\in
N(H,u)$. Then:
\begin{enumerate}
\item[{\rm 1.}]
$\gamma\in N_G^1(H^0)$;
\item[{\rm 2.}]
$\pi(H)=\pi(\gamma H\gamma\inv)$;
\item[{\rm 3.}]
$g\gamma\in N^\ast(H,u)$; and
\item[{\rm 4.}]
$N(H,u)$ contains an open closed subset of $N(\gamma H \gamma\inv,u)$
containing $g$.
\end{enumerate}
\end{proposition}

\proof Replacing $u$ by $g\inv u g$, we may assume that
$g=\{e\}$. Since  $e\in N^\ast(H,u)$ and $\gamma\in N(H,u)$, by
Proposition~\ref{nimish:prop:g:in:N*(H,u)},
\[
\pi(H)=\cl{\pi(\langle u\rangle)}=\cl{\pi(\langle u\rangle\gamma)}
\subset \gamma\pi(H)=\pi(\gamma H\gamma\inv).
\]
By the dimension consideration, $\pi(H^0)=\pi(\gamma H^0\gamma\inv)$, and hence
$H^0=\gamma H^0\gamma\inv$. Since the action of $\gamma$ on $G/\Gamma$
is a homeomorphism, the number of connected components of $\pi(H)$ and
$\gamma\pi(H)$ are the same. Therefore $\pi(H)=\gamma\pi(H)$. Hence
$\gamma\in N^\ast(H,u)$. Moreover, since with respect to the Haar
measures, $\vol(\pi(\gamma H^0 \gamma\inv))=\det(\Ad_G(\gamma)|_{{\rm
Lie}(H^0)})\vol(\pi(H^0))$, we obtain that $\gamma\in
N^1_G(H^0)$. 

By definition, the sets $$\{h\in G^0: h\inv u h\in uH^0\} \mbox{~and~} \{h\in G^0:
h\inv u h\in u\gamma H^0\gamma\inv\}$$ are open closed in
$N(H,u)$ and $N(\gamma H\gamma\inv,u)$, respectively. Since  $\gamma\in  
N_G^1(H^0)$, statement~(4) follows. \qed

\begin{proposition}  \label{nimish:prop:measure:on:tube}
Let $H\in\cH_u$ and $\lambda$ be a $u$-invariant $u$-ergodic Borel
probability measure on $\pi(N^\ast(H,u))$. Then $\lambda=g\lambda_H$,
for any $$g\in N^\ast(H,u)\cap \pi\inv(\supp \lambda).$$
\end{proposition}

\proof
By Hedlund's lemma, there exists $g_0\in N^\ast(H,u)$ such that $$\cl{\langle
u\rangle\pi(g_0)}=\supp(\lambda).$$ Therefore by
Proposition~\ref{nimish:prop:g:in:N*(H,u)}, $\supp(\lambda)=g_0\pi(H)$. Now by
Ratner's theorem as discussed in the preceding subsection, we have
that $\lambda$ is $g_0Hg_0\inv$-invariant. Hence
$\lambda=g_0\lambda_H$.

Now $\pi\inv(\supp\mu)=g_0H\Gamma$. By Proposition~\ref{nimish:prop:N*(H,u)},
if $g\in g_0H\Gamma\cap N^\ast(H,u)$, then $g\pi(H)=g_0\pi(H)$. Hence
$g\lambda_H=g_0\lambda_H=\lambda$. This completes the proof of the
proposition.  \qed

\section{Abundance of Unipotent Subgroups}

The following main technical result of this section is used in the proofs 
of Theorems~\ref{nimish:thm:measure} and \ref{nimish:thm:closure}. 

\begin{proposition}  \label{nimish:prop:unip}
Let $G$ be a Lie group, $H$ a closed connected subgroup of $G$, $u\in
N_G(H)$ an $\Ad_G$-unipotent element, and $U$ the subgroup generated
by all one-parameter $\Ad_G$-unipotent subgroups of $H$. Then the set
$$\{g\in N_G(U): g\inv ugu\inv\in U\}$$ contains a neighbourhood of $e$
in the set $$\{g\in N_G(U):g\inv ugu\inv\in H\}.$$
\end{proposition}

\proof Let $\tilde{U}=\Ad_G(U)$. Then $\tilde{U}$ is a connected real
algebraic group generated by one-parameter unipotent subgroups of
$\GL(\lca{G})$ (cf.~\cite[Proof of Lemma~2.9]{Shah:distribution}). Put
$\tilde{L}=N_{\GL(\lca{G})}(\tilde{U})$. Then $\tilde{L}$ and
$\tilde{L}/\tilde{U}$ are real algebraic groups, and the natural
quotient homomorphism $\tilde{q}:\tilde{L}\to\tilde{L}/\tilde{U}$ is
algebraic. Put $L=N_G(U)$. 

\begin{claim} \label{claim:U}
There exists a neighbourhood $\Omega$ of $e$ in $L$ such that for any
$h\in H\cap\Omega$, if $\tilde q(\Ad_G(h))$ is an algebraic unipotent element
of $\tilde L/\tilde U$, then $h\in U$. 
\end{claim}

To show this, let $\tilde\Omega$ be a neighbourhood of the identity in
$\tilde{L}/\tilde{U}$ such that the following holds: For any
one-parameter subgroup ${\bar{h}(t)}\subset \tilde{L}/\tilde{U}$, if
$\bar{h}([0,1])\subset\tilde\Omega$ and $\bar{h}(1)$ is an algebraic
unipotent element, then $\{\bar{h}(t)\}$ is algebraic unipotent
subgroup. In this case, there exists a unipotent one-parameter
subgroup $\{h(t)\}\subset\tilde{L}$ such that
$\tilde{q}(h(t))=\bar{h}(t)$. Note that 
$$
H \subset L=N_G(U)\subset \Ad_G\inv(\tilde L).
$$ 
Let $\Omega_1=\Ad^{-1}_G \circ \tilde q\inv(\tilde\Omega)$. Let
$\Omega\subset\Omega_1$ be a neighbourhood of $e$ in $L$ such that the
following holds: for any $h\in H\cap\Omega$, there exists a
one-parameter subgroup $\{h(t)\}\subset H$ such that $h=h(1)$ and
$h([0,1])\subset\Omega_1$. Now since $\tilde U\subset\Ad_G(H)$, the
claim follows from the above construction. \qedclaim

\bigskip
Let $\lca{L}$ and $\lca{U}$ denote the Lie algebras
corresponding to $L$ and $U$, respectively. We identify the Lie
algebra of $L/U$ with $\lca{L}/\lca{U}$. Let $q:L\to L/U$ be the natural
quotient homomorphism.

Now suppose that the proposition is not true. Then there exists a
sequence $g_k\to e$ in $L$ such that $g_k\inv u g_k u\inv\in
H\setminus U$ for all $k\in\NN$.  By passing to a subsequence, there
exists $X_k\in\lca{L}/\lca{U}$ such that $$q(g_k)=\exp_{L/U}(X_k)$$ for
all $k\in\NN$, and $X_k\to 0$ as $k\to\infty$.

Consider the linear action of $L$ on $\lca{L}/\lca{U}$ via the
representation $\Ad_{L/U}\circ q$. Since $u g_ku\inv\not\in g_kU$, we
have that $u\cdot X_k\neq X_k$. Now since $\Ad_{L/U}(q(u))$ is
unipotent, $u^n\cdot X_k\to\infty$ as $n\to\infty$. Therefore there
exists a sequence $n_k\to\infty$ such that after passing to a
subsequence $u^{n_k}\cdot X_k\to X$, where $X\in\lca{L}/\lca{U}\setminus
0$. We can choose $\{n_k\}$ such that $\exp_{L/U}(X)=q(h)\neq e$
for some $h\in\Omega$. Thus $q(u^{n_k}g_ku^{-n_k})\to q(h)$ as
$k\to\infty$. Now $U\subset H$, $g_k\inv ug_ku\inv\in H$ and
$g_k\to e$ as $k\to\infty$. Therefore $h\in H\setminus U$.

Now 
\[
(\tilde{q}(\Ad_G u))^{n_k} \tilde q(\Ad_G g_k) (\tilde q(\Ad_G
u))^{-n_k}\to \tilde{q}(\Ad_G h) \qquad \mbox{as $k\to\infty$}.
\]
Since $\tilde{q}(\Ad_G g_k)\to e$ as $k\to\infty$, it follows that
$\tilde{q}(\Ad_G h)$ is an algebraic unipotent
element of $\tilde{L}/\tilde{U}$. Hence by Claim~\ref{claim:U},
$h\in U$, which is a contradiction.
\qed

The following simple observation is useful. 

\begin{lemma} \label{lemma:analytic}
Let $G$ be a Lie group (so that $G^0$ is analytic), $F$ a  connected
(analytic) Lie group, and $\rho:F\to G^0$ an analytic map. Let
$\sigma_F$ denote a Haar measure on $F$, and $M$ be an analytic
submanifold of $G^0$.  If $\sigma_F(\rho\inv(M))>0$, then
$\rho(F)\subset M$.

In particular, if $F$ is an analytic subgroup of $G$ such that
$N(H,u)$  contains a neighbourhood of $e$ in $F$ then
$F\subset N(H,u)$, where $N(H,u)$ is as defined in 
Section~\ref{sec:singular:subsets}.
\end{lemma}

\begin{corollary} \label{nimish:cor:normal}
Let $G$ be a Lie group, $H$ a closed connected subgroup of $G$, $u\in
N_G(H)$ an $\Ad_G$-unipotent element, and $U$ be a the closed connected
subgroup of $H$ generated by all one-parameter $\Ad_G$-unipotent
subgroups of $H$. Then $h\inv uhu\inv \in U$ for all $h\in H$. In
particular, $u\in N_G(F)$ for any subgroup $F$ of $H$ containing $U$.
\end{corollary}

\proof Let $\rho: H \to G$ be the map defined by $$\rho(h)=h\inv u h
u\inv \mbox{~for~ all~} h\in H.$$ Since $\rho(H)\subset H$, by
Proposition~\ref{nimish:prop:unip}, $\rho\inv(U)$ contains a neighbourhood of
$e$ in $H$. Therefore by Lemma~\ref{lemma:analytic}, $h\inv u h u\inv
\in U$ for all $h\in H$. \qed

\section{Proofs of Theorem~\ref{nimish:thm:measure} and
Theorem~\ref{nimish:thm:locally:finite}}

For locally finite $u$-invariant measures, we have the
following result due to Dani:

\newcommand{\citone}{\cite[Theorem 4.3]{Dani:rk=1}}

\begin{theorem}[\mbox{\citone}] \label{nimish:thm:property(D)} 
Let $G$ be a Lie group and $\Gamma$ a closed subgroup such that
$G/\Gamma$ has a finite $G$-invariant measure. Let $u\in G$ be an
$\Ad_G$-unipotent element, and $\mu$ be a $u$-invariant locally finite
Borel  measure on $G/\Gamma$. Then there exist a partition of
$G/\Gamma$ into countably many $u$-invariant Borel measurable subsets,
say $X_i$ $(i\in\NN)$, such that $\mu(X_i)<\infty$, $\forall \ i\in\NN$.
\end{theorem}

Due to this result, the Theorem~\ref{nimish:thm:measure} and
Theorem~\ref{nimish:thm:locally:finite} are special cases of the following.

\begin{theorem} \label{nimish:thm:measure:locally:finite}
Let $G$, $\Gamma$, $W$, and $\cU$ be as in
Theorem~\ref{nimish:thm:measure}. Let $\mu$ be a locally finite $W$-invariant
$W$-ergodic Borel measure on $G/\Gamma$. Suppose that for any
$\Ad$-unipotent element $u\in\cU$, there exists a partition of
$G/\Gamma$ into countably many Borel measurable $u$-invariant subsets
$X_i$ $(i\in\NN)$ such that $\mu(X_i)<\infty$, $\forall i\in\NN$. Then
there exists a closed subgroup $H$ of $G$ containing $W$ such that
$\supp(\mu)$ is closed $H$-orbit.
\end{theorem}

We intend to prove this theorem by induction on the dimension of
$G^0$. Note that the theorem is obvious, if $G$ is a discrete group;
that is, if $\dim(G^0)=0$.

The rest of the proof of this theorem is a series of claims and
propositions. 

\begin{claim} \label{claim:zero} 
We may assume that $\cU$ is finite. 
\end{claim}

\proof Since the Zariski closure of a cyclic group generated by a
unipotent linear transformation is a connected group, we have that the
Zariski closure of $\langle \Ad_G\cU\rangle$ is a connected real
algebraic group of dimension, say $d$. Hence there is a subset
$\cU_1\subset \cU$ consisting of at most $d$ elements such that
$\Zcl(\langle \Ad_G\cU_1\rangle)=\Zcl(\langle \Ad_G \cU
\rangle)$. Thus without loss of generality\ we may replace $\cU$ by $\cU_1$ and assume
that $\cU$ is finite. \qed

Let $\pi:G\to G/\Gamma$ denote the natural quotient map. 

\begin{claim} \label{claims:e}
We may assume that for each $u\in\cU$, there exists a $u$-invariant
Borel measurable subset $X_u$ of $G/\Gamma$ such that
$\mu(X_u)<\infty$, and $\pi(e)$ belongs to the support of the
restriction of $\mu$ to $X_u$, and $\supp(\mu) \subset \cl{\pi(W)}$. 
\end{claim}

\proof Using Hedlund's lemma and Claim \ref{claim:zero} it is straightforward
to  obtain the conclusion of the claim for
$\pi(g)$ in place of $\pi(e)$ for some $g$ in $G$. Now if we work with
$g\Gamma g\inv$, in place of $\Gamma$, without loss of generality\ we may assume\ that $\pi(e)$
belongs to the support of the restriction of $\mu$ to $X_u$.  \qed

\begin{claim} \label{claims:two} 
We may assume that $\Gamma$ is a discrete subgroup of $G$.
\end{claim}

\proof For each $u\in\cU$, by Claim~\ref{claims:e}, $\pi(e)$ belongs
to the support of a finite $u$-invariant Borel measure on $G/\Gamma$. 
Therefore by Proposition~\ref{nimish:prop:closed:discrete}, we
have that $\Ad_G(\langle u \rangle)\subset \Zcl(\Ad_G(\Gamma))$ for
all $u\in\cU$. Therefore $W\subset N_G(\Gamma^0)$. Since $\supp(\mu) \subset 
\cl{\pi(W)}$, replacing
$G$ by $N_G(\Gamma^0)$, without loss of generality\ we may assume\ that $\Gamma^0$ is normal in
$G$. Now again replacing $G$ by $G/\Gamma^0$ and $\Gamma$ by
$\Gamma/\Gamma^0$, the claim holds.  \qed

\begin{claim} \label{claims:minimal}
We may assume that if $\mu(\pi(F))>0$ for any closed connected
subgroup $F$ of $G$ such that $W\subset N_G(F)$, then
$F=G^0$.
\end{claim}

\proof For any $w_1,w_2\in W$, either $$w_1\pi(F)=w_2\pi(F) \mbox{~or~}
w_1\pi(F)\cap w_2\pi(F)=\emptyset.$$ Since $\mu$ is locally finite,
and $\mu(\pi(F))>0$, the group $WF/F$ is countable. Put $H=WF$, and
extend the topology of $F$ to $H$ such that $H^0=F^0$. Let
$\Lambda=H\cap\Gamma$, and consider the natural continuous inclusion
$\rho:H/\Lambda\to G/\Gamma$. By the ergodicity of $W$-action, $\mu$
is concentrated on $\pi(H)$. Hence $\mu$ can be treated as a locally
finite $W$-invariant $W$-ergodic Borel probability measure on
$H/\Lambda$. Suppose that $F\neq G^0$. Then $\dim(H^0)<\dim(G^0)$, and
hence by induction hypothesis there exists a closed subgroup $H_1$ of
$H$ containing $W$ such that $\mu$ is $H_1$-invariant and $\supp(\mu)$
is a closed $H_1$-orbit in $H/\Lambda$. Thus $H_1$ is a measurable
subgroup of $G$ containing $W$ such that $\mu$ is $H_1$-invariant 
and concentrated on an
orbit of $H_1$. Now the theorem follows from
Lemma~\ref{lemma:finitevol:closed}. Thus without loss of generality\ we may assume\ that
$F=G^0$. \qed

\begin{claim} \label{claims:one} 
We may assume that $G=WG^0=G^0\Gamma$.
\end{claim}

\proof 
Since $\pi(WG^0)=W\pi(G^0)$ is a closed $W$-invariant subset of
$G/\Gamma$ with strictly positive $\mu$-measure, by ergodicity, we have
that $$\supp(\mu)\subset\pi(WG^0).$$ Therefore
replacing $G$ by $WG^0$, without loss of generality\ we may assume\ that $G=WG^0$.

For any $w_1,w_2\in W$, either $$w_1\pi(G^0)=w_2\pi(G^0) \mbox{~or~}
w_1\pi(G^0)\cap w_2\pi(G^0)=\emptyset.$$ Let $\mu_0$ be the
restriction of $\mu$ to $\pi(G^0)$. Define 
$$
W_0=\{w\in W:w\mu_0=\mu_0\} = \{w \in W : \pi(w) \in \pi(G^0)\}. 
$$ 
Let $u\in\cU$. By Claim~\ref{claims:e},
$\mu_0(X_u)>0$. Therefore there exists $k\in\NN$ such that $u^k\in
W_0$. Note that $\Zcl(\langle \Ad_G u \rangle)=\Zcl(\langle \Ad_G u^k
\rangle)$. Put
\[
\cU_0=W_0 \bigcap \left( \bigcup_{u\in\cU} \langle u \rangle \right).
\]
Then $\cU_0$ consists of $\Ad_G$-unipotent elements, and 
\begin{equation} \label{eq:U_0}
\Ad_G(W)\subset \Zcl(\langle \Ad_G \cU_0 \rangle).
\end{equation} 

Now suppose we can prove the theorem for the action of $W_0$ on
$\mu_0$. Then there exists a closed subgroup, say $F$, such that
$W_0\subset F$, $\mu_0$ is $F$-invariant and $\supp(\mu_0)=\pi(F)$. By
Equation~\ref{eq:U_0}, we have that $W\subset N_G(F^0)$. Therefore by
Claim~\ref{claims:minimal}, $F^0=G^0$. Now since $w\mu_0$ is
$G^0$-invariant for any $w\in W$, we have that $\mu$ is
$G^0$-invariant. Thus $\mu$ is $G$-invariant, and the theorem
follows. Thus without loss of generality\ we may assume\ that $G/\Gamma=\pi(G^0)$. \qed

\begin{proposition} \label{nimish:prop:capN}
There exist $H_u\in\cH_u$ for all $u\in\cU$ such
that $$\mu(\pi(\capsu))>0.$$
\end{proposition}

\proof Let $\Omega\subset G$ be a Borel measurable subset such that
$\mu(\pi(\Omega))>0$. Express $\cU=\{u_1,\ldots,u_k\}$ for some
$k\in\NN$. 
There exists a $u_1$-invariant Borel
measurable subset $X_1$ of $G/\Gamma$ such that $\mu(X_1)<\infty$, and
$\mu(\pi(\Omega)\cap X_1)>0$. Let $\mu_1$ denote the restriction of
$\mu$ to $X_1$. Then $\mu_1$ is a finite $u_1$-invariant measure.

Now $\mu_1$ is a direct integral of finite $u_1$-ergodic
$u_1$-invariant measures; (see \cite[Section~1.4]{Dani:invent78} for a
precise statement). Therefore by Proposition~\ref{nimish:prop:measure:tube}
and Proposition~\ref{nimish:prop:countable}, there exists $H_1\in\cH_{u_1}$
such that $$\mu_1(\pi(N^\ast(H_1,u_1))\cap\pi(\Omega))>0.$$ Therefore
there exists $\gamma_1\in\Gamma$ such that
\[
\mu_1(\pi(N^\ast(H_1,u_1)\gamma_1\cap\Omega))>0.
\]
Put $H_{u_1}=\gamma_1\inv H_1\gamma_1$ and
$\Omega_1=N^\ast(H_{u_1},u_1)\cap\Omega$. Then $\mu(\pi(\Omega_1))>0$.

For $2\leq i\leq k$, we inductively carry out the same procedure with
$u_i$ in place of $u_1$, and $\Omega_{i-1}$ in place of $\Omega$. We
obtain $H_{u_i}\in\cH_{u_i}$ such that if we put
$\Omega_i=N^\ast(H_{u_i},u_i)\cap\Omega_{i-1}$, then
$\mu(\pi(\Omega_i))>0$.

Now $\Omega_k\subset\capsu$ and $\mu(\pi(\Omega_k))>0$. This
completes the proof.
\qed

\begin{proposition} \label{nimish:prop:analytic}
Let $F$ be a connected Lie subgroup of $G$ such that $\pi(F)$ has
an $F$-invariant probability measure, say $\lambda_F$. Suppose for
some $g\in G$, $(g\lambda_F)(\pi(\capu))>0$. Then there exists
$\gamma\in\Gamma$ such that
\[
g F\gamma\subset\capu.
\]
\end{proposition}

\proof
Let $\tilde{\lambda}_F$ denote a Haar measure on $F$. Then 
\[
\tilde\lambda_F(\{h\in F:g\pi(h)\in\pi(\capu)\})>0.
\]
Therefore there exists $\gamma\in\Gamma$ such that 
\[
\tilde\lambda_F(\{h\in F:gh\gamma\in\capu\})>0.
\]

Take any $u\in\cU$. Then the map $\rho:F\to G$, given by
$$\rho(h)=(gh\gamma)\inv u (gh\gamma) u^{-1} \mbox{~for all~} h\in F,$$ is
analytic. Since $\tilde\lambda(\rho\inv(H_u))>0$, by
Lemma~\ref{lemma:analytic} $\rho(F)\subset H_u$. Hence
$gF\gamma\subset N(H_u,u)$. This completes the proof.  \qed

\begin{proposition}   \label{nimish:prop:+:in}
For any $g\in\capsu$ and $v\in\cU$, if
$$g\lambda_{H_v^0}(\pi(\capu))>0$$ then
\begin{equation}  \label{eq:in}
gH_v^0\subset\capu
\end{equation}
and
\begin{equation}  \label{eq:1}
g\lambda_{H_v^0}(\pi(\cup_{u\in\cU} S(H_u,u)))=0.
\end{equation}
\end{proposition}

\proof
By Proposition~\ref{nimish:prop:analytic}, there exists $\gamma\in\Gamma$ such
that 
\[
gH_v^0\gamma\subset\capu.
\]
Hence 
\[
gH_v^0\subset\cap_{u\in\cU}N(\gamma H_u \gamma\inv,u).
\]

By Proposition~\ref{nimish:prop:N*(H,u)}, the set $N(H_u,u)$ contains the
connected component of $g$ in $N(\gamma
H_u\gamma\inv,u)$. Therefore $gH_v^0\subset N(H_u,u)$ for all
$u\in\cU$. Hence~(\ref{eq:in}) holds.

Now suppose $g\lambda_{H_v^0}(\cup_{u\in\cU}\pi(S(H_u,u)))>0$. Then
there exist $u\in\cU$ and $F_u<H_u$ such that
$g\lambda_{H_v^0}(\pi(N(F_u,u)))>0$. By
Proposition~\ref{nimish:prop:analytic}, there exists $\gamma\in\Gamma$ such
that $gH^0_u\gamma\subset N(F_u,u)$. Therefore $g\gamma\in
N(F_u,u)$. Thus $g\gamma\not\in\Nsu$. Since $g\in\Nsu$, this
contradicts Proposition~\ref{nimish:prop:N*(H,u)}.  \qed

\begin{proposition}  \label{nimish:prop:product}
There exists a Borel set $T^\ast\subset\capsu$ such that
\[
\mu(\pi(\capsu\setminus T^\ast))=0 \mbox{~and~}
g\lambda_{H_u^0}(\pi(T^\ast))=1,\ \forall u\in\cU,\ \forall g\in T^\ast.
\]
\end{proposition}

\proof Put $T^\ast_0=\capsu$. Let $\mu^\ast$ be the restriction of
$\mu$ to $\pi(T^\ast_0)$. Let
\[
T^\ast_1=\{g\in T^\ast_0:g\lambda_{H^0_u}(\pi(T^\ast_0))>0,\,\forall
u\in\cU\}.
\]
Considering the ergodic decomposition of $\mu^*$ restricted to 
$\pi(N(H_u,u))$ for all $u\in \cU$, and applying Proposition~\ref{nimish:prop:measure:on:tube}, we
obtain that $\mu(\pi(T^\ast_0\setminus T^\ast_1))=0$. By
Proposition~\ref{nimish:prop:+:in},
\[
T^\ast_1=\{g\in T^\ast_0:g\lambda_{H^0_u}(\pi(T^\ast_0))=1,\,\forall
u\in\cU\}. 
\]
Therefore $\mu$ restricted to $\pi(T_1^\ast)$ is $\mu^\ast$. Now for
any $u\in\cU$, by Proposition~\ref{nimish:prop:measure:on:tube}, the measure
$\mu^\ast$ is a direct integral of measures of the form
$g\lambda_{H^0_u}$, where $g\in T_1^\ast$.

By the above procedure we can obtain a decreasing sequence of Borel
subsets $\{T^\ast_i\}$ such that 
\[
T^\ast_{i+1}=\{g\in
T^\ast_i:g\lambda_{H^0_u}(\pi(T^\ast_i))=1,\,\forall u\in\cU\},
\]
and $\mu(\pi(T^\ast_i\setminus T^\ast_{i+1}))=0$ for all $i\geq 0$. 
Thus $T^\ast=\cap_{i\geq 0} T^\ast_i$ has the desired properties.
\qed

\begin{claim} \label{claims:generic}
We may assume that $e\in\capsu$, $\cl{\pi(W)}=\supp(\mu)$, and
$\pi(e)$ is in the support of the restriction of $\mu$ to $\pi(\capsu)$.  
\end{claim}

\proof Let $T^\ast$ be as in Proposition~\ref{nimish:prop:product}. By
Hedlund's lemma, there exists $g\in T^\ast$ such that
$\cl{W\pi(g)}=\supp(\mu)$ and $\pi(g)$ is a density point of
$\mu^\ast$. Replacing $\Gamma$ by $g \Gamma g\inv$, and $H_u$ by
$gH_ug\inv$ for all $u\in\cU$, without loss of generality\ we may assume\ that $e\in T^\ast$,
$\cl{\pi(W)}=\supp(\mu)$, and $\pi(e)$ is a density point of
$\mu^\ast$. 
\qed

\begin{proposition} \label{nimish:prop:Hv} 
$H_v^0\subset\capu$ for all $v\in\cU$.
\end{proposition}

\proof Since $e\in T^\ast$, by Proposition~\ref{nimish:prop:product}, we have
that $$\lambda_{H_v^0}(\pi(\capu))=1$$ for all $v\in\cU$. Now the
proposition follows from Proposition~\ref{nimish:prop:+:in}.
\qed

\begin{definition} \label{not:Uu} \rm
For each $u\in\cU$, let $U_u$ be the subgroup generated by all
one-parameter $\Ad_G$-unipotent subgroups of $H_u$. Then $U_u$ is
normal in $H_u$. Let $F_u$ be the connected component of the identity
in the closure of the subgroup $U_u(H_u\cap\Gamma)$. Note that
$\cl{\pi(U_u)}=\pi(F_u)$.
\end{definition}

\begin{proposition} \label{nimish:prop:normalize}
\[
F_v\,\subset\, N_G(H_u^0) \cap N_G(F_u)\,\subset\, N_G(U_u).
\]
\end{proposition}

\proof Since $e\in \capsu$, by Proposition~\ref{nimish:prop:N*(H,u)},
$$\Gamma\cap N(H_u,u)\subset N_G(H_u^0).$$ Note that $\Gamma\cap
N_G(U_u)\subset N_G(F_u)$. Now since $N_G(H_u^0)\subset N_G(U_u)$, we
have $$\Gamma\cap N(H_u,u)\subset N_G(F_u).$$ Let $v\in\cU$. By
Proposition~\ref{nimish:prop:Hv}, $F_v\subset H^0_v\subset\capu$. Hence
$F_v\cap\Gamma\subset N_G(H_u^0)\cap N_G(F_u)$. By
Proposition~\ref{nimish:prop:closed:discrete},
$\Zcl(\Ad_G(F_v\cap\Gamma))\supset\Ad_G(F_v)$. Now the proposition
follows.  \qed

\begin{proposition} \label{nimish:prop:F}
Express $\cU=\{u_1,\ldots,u_k\}$ and put $$F=\{f_1\cdots f_k\in
G:f_i\in F_{u_i},\, 1\leq i\leq k\}.$$ Then $F$ is a closed subgroup of
$G$, $\pi(F)$ is closed, and $F\cap\Gamma$ is a lattice in $F$.
\end{proposition}

\proof Put $F_i=F_{u_i}$. Due to Proposition~\ref{nimish:prop:normalize}, we
can define the semidirect product $\tilde F=F_k \times\cdots\times
F_1$, where $F_i$ acts on $F_{i-1}\times\cdots\times F_1$ by
conjugation on each factor ($i=2,\ldots,k$). Clearly $\tilde F$ is a
connected Lie group. Let $\Lambda=(F_k\cap\Gamma)\times\dots\times
(F_1\cap\Gamma)$. Since $F_i\cap\Gamma$ is a lattice in $F_i$, we have
that $\Lambda$ is a lattice in $\tilde F$. Let $\tilde\sigma$ be a
finite $\tilde F$-invariant measure on $\tilde F/\Lambda$. By the
definition of semidirect product, the map $\rho:\tilde F\to G$, given
by $\rho(f_k,\ldots,f_1)=f_k\cdots f_1$ for all $f_i\in F_i$
($i=1,\ldots,k$), is a continuous  homomorphism. Note that $F=\rho(\tilde
F)$.  Since $\rho(\Lambda)\subset \Gamma$, the map $\rho$  determines a
continuous $\rho$-equivariant map $\bar\rho:\tilde F/\Lambda\to
G/\Gamma$. Then the push-forward of $\tilde\sigma$ under $\bar\rho$ is
a finite $F$-invariant measure concentrated on $\pi(F)=\bar\rho(\tilde
F/\Lambda)$. Now by Lemma~\ref{lemma:finitevol:closed}, $F$ is closed
and $\pi(F)$ is closed. Since $\pi(F)$ has a finite $F$-invariant
measure, $F\cap\Gamma$ is a lattice in $F$.  \qed

\begin{proposition}  \label{nimish:prop:W:in:N(F)}
$W\subset N_G^1(F)$.
\end{proposition}

\proof Take any $u,v\in\cU$. Define the map $\rho:F_v\to G$ by
$$\rho(f)=f\inv u f u\inv \mbox{~for~ all~} f\in F_v.$$ Since $F_v\subset
N(H_u,u)$ and $F_v$ is connected, we have $\rho(f)\in H_u^0$ for all
$f\in F_v$. Now since $F_v\subset N_G(U_u)$, by
Proposition~\ref{nimish:prop:unip}, $\rho\inv(U_u)$ contains a neighbourhood
of $e$ in $F_v$. Therefore by Lemma~\ref{lemma:analytic}, $f\inv
ufu\inv\in U_u\subset F$, $\forall f\in F_v$. Hence $u\in N_G(F)$. Now
the proposition follows.  \qed

\begin{claim} \label{claims:F:normal}
We may assume that $F$ is a normal subgroup of $G$.
\end{claim}

\proof
By Proposition~\ref{nimish:prop:nor}, $\pi(N_G^1(F))$ is closed. Now by
Proposition~\ref{nimish:prop:W:in:N(F)}, we have 
$\supp(\mu)=\cl{\pi(W)}\subset\pi(N_G^1(F))$. Therefore without loss of generality\ we
can replace $G$ by $N_G^1(F)$. Now the claim follows.
\qed

\begin{proposition}   \label{nimish:prop:F.mu=mu}
The measure $\mu$ is $F$-invariant.
\end{proposition}

\proof For any $u\in\cU$, since $F_u\subset H_u^0$ and $\pi(F_u)$
has a finite $F_u$-invariant measure, by
Proposition~\ref{nimish:prop:measure:on:tube} and
Proposition~\ref{nimish:prop:product}, $\mu^\ast$ is a direct integral of
measures of the form $g\cdot\lambda_{F_u}$, where $g\in T^\ast$.
Therefore by arguing as in Proposition~\ref{nimish:prop:product} for $F_u$ in
place of $H_u^0$, we obtain a Borel set $T'\subset T^\ast$ such that
the following holds: $\mu^\ast(\pi(T^\ast\setminus T'))=0$, and
$g\lambda_{F_u}(\pi(T'))=1$, $\forall g\in T',\, \forall
u\in\cU$. Hence $g\cdot\lambda_F(\pi(T'))=1$ $(\forall g\in T'$), and
the measure $\mu^\ast$ is a direct integral of measures of the form
$g\cdot\lambda_F$ ($g\in T'$).

Since $F$ is a normal subgroup of $G$, the measure $g\lambda_F$ is
$F$-invariant for all $g\in G$. Therefore $\mu^\ast$ is
$F$-invariant. Hence $w\mu^\ast$ is $F$-invariant for all $w\in
W$. Since $\mu$ is $W$-ergodic, and $\mu^\ast\neq 0$, we conclude that
$\mu$ is $F$-invariant. 
\qed

\begin{claim}  \label{claims:mu(Z_G(W))}
We may assume  that $\mu(\pi(Z_G(W)^0))>0$.
\end{claim}

\proof Due to Claim \ref{claims:F:normal} and Proposition~\ref{nimish:prop:F.mu=mu}, without loss of generality\ we can pass to
the quotient of $G$ by $F$ and assume that $F=\{e\}$. 

Take any $u\in\cU$. Since $U_u=\{e\}$, by Proposition~\ref{nimish:prop:unip},
there exists a neighbourhood $\Omega_u$ of $e$ in $N(H_u,u)$ such that
$\omega\inv u\omega u\inv\in\{e\}$ for all $\omega\in\Omega_u$. Thus
$Z_G(u)^0$ contains a neighbourhood of $e$ in $N(H_u,u)$. Therefore
$Z_{G}(\langle\cU\rangle)^0$ contains a neighbourhood of $e$ in
$\capu$. Now since $Z_G(W)^0=Z_G(\langle \cU \rangle)^0$, by
Claim~\ref{claims:generic}, we have that $\mu(\pi(Z_G(W)^0))>0$.  \qed

\begin{claim} \label{claims:three}
$G^0\subset Z_G(W)$.
\end{claim}

\proof
This follows from Claims~\ref{claims:minimal} and
\ref{claims:mu(Z_G(W))}.
\qed

\begin{claim}
We may assume that $G^0\cap\Gamma$ is contained in the center of $G$.
\end{claim}

\proof By Lemma~\ref{lemma:centralizer}, the orbit
$\pi(Z_G(\gamma))$ is closed for any $\gamma\in\Gamma$. By
Claim~\ref{claims:three}, for any $\gamma\in G^0\cap \Gamma$, we have
$W\subset Z_G(\gamma)$. Therefore
$$\supp(\mu)=\cl{\pi(W)}\subset\pi(Z_G(\gamma)).$$ By
Claim~\ref{claims:minimal}, $G^0\subset Z_G(\gamma)$, and by
Claim~\ref{claims:one} $Z_G(\gamma)=G$.  \qed

\bigskip

\noindent{\bf Completion of the proof of the theorem} 
Put $Z=G^0/G^0\cap\Gamma$. Then $Z$ is a locally compact
group. Consider the natural inclusion $\psi:Z\to G/\Gamma$. Since
$G=G^0\Gamma$, the map $\psi$ is a homeomorphism. Let the map $\rho:W\to
Z$ be defined by $\rho(w)=\psi\inv(\pi(w))$ ($\forall w\in W$). 
By Claim~\ref{claims:three}, we have $w\psi(z)=\psi(z
\rho(w))$ for all $z\in Z$. Therefore the action of $w\in W$ on $G/\Gamma$
corresponds to the right action of $\rho(w)$ on $Z$. Let $\tilde{\mu}$
be the projection of $\mu$ under $\psi\inv$. Then $\tilde{\mu}$ is
ergodic under the right action of $\rho(W)$ on the locally compact
group $Z$. By Claim~\ref{claims:generic}, $\supp(\tilde{\mu})$
contains the identity element of $Z$. Hence $\tilde{\mu}$ is a Haar
measure on the closed subgroup $\cl{\rho(W)}$. Let $H$ be the inverse
image of $\cl{\rho(W)}$ in $G^0$. Since $\mu=\psi_\ast(\tilde{\mu})$,
we have that $\mu$ is $H$-invariant and $\supp(\mu)=\pi(H)$. This
completes the proof of Theorem~\ref{nimish:thm:measure:locally:finite}.  \qed

\section{Orbit Closures in Finite Volume Homogeneous Spaces}

One of the main purposes of this section is to prove the following
result, which will be used in the proof of
Theorem~\ref{nimish:thm:closure}. 

\begin{theorem} \label{nimish:thm:closed:finite}
Let $G$, $\Gamma$, and $W$ be as in Theorem~\ref{nimish:thm:closure}. Let
$\pi:G\to G/\Gamma$ be the natural quotient map and $x_0=\pi(e)$. Let
$H$ be the minimal among the closed subgroups $F$ of $G$ such that
$W\subset F$ and $Fx_0$ is closed. Then $H^0x_0$ has a finite
$H^0$-invariant measure.
\end{theorem}

To prove this theorem we will reduce this question to homogeneous
spaces of semisimple groups. In order to quotient out by solvable
factors, first we study the effects on the orbit closures when we pass
to the quotients of finite volume homogeneous spaces.

\begin{lemma} \label{lemma:central:project}
Let $G$ be a locally compact group, $\Lambda$ a closed subgroup $G$,
and put $x_0=e\Lambda$. Let $Z$ be a closed subgroup of $\{g\in
G:g\Lambda g\inv\subset\Lambda\}$ 
such that $Zx_0$ is compact. Let $\rho:G/\Lambda \to
G/Z\Lambda$ be the natural quotient map. Suppose there exists a closed
subgroup $H$ of $G$ such that the orbit $Hx_0$ is open in its closure
$\cl{Hx_0}$. Then
\[
\rho(\cl{Hx_0}\setminus Hx_0)=\cl{H\rho(x_0)}\setminus H\rho(x_0).
\]
In particular, if $H\rho(x_0)$ is closed, then $Hx_0$ is closed.  
\end{lemma}

\proof
Since $Z\Lambda/\Lambda$ is compact, $\rho$ is a proper map. Therefore
$$\rho(\cl{Hx_0}\setminus Hx_0)\supset\cl{H\rho(x_0)}\setminus H\rho(x_0).$$

To show the inclusion, suppose that there exists $y_0\in
Y=\cl{Hx_0}\setminus Hx_0$ such that $\rho(y_0)=h\rho(x_0)$ for some
$h\in H$. Then there exists $z\in Z$ such that $y_0=hz\Lambda$. Hence
$zx_0\in Y$. 

We claim that $z^kx_0\in Y$ for all $k\in\NN$. To prove this claim by
induction, suppose that $z^kx_0\in Y$. Now sequences $\{h_i\}\subset
H$ and $\lambda_i\in\Lambda$ be such that $h_i\lambda_i\to z$ as
$i\to\infty$. Then 
\[
h_iz^kx_0=(h_i\lambda_i)\lambda_i\inv z^kx_0 = (h_i\lambda_i)z^kx_0 
\to zz^kx_0,
\]
as $i\to\infty$. Since $Y$ is a closed $H$-invariant set, we have that
$z^{k+1}x_0\in Y$. This proves the claim. 

Since $Zx_0$ is compact, and $Z\cap\Lambda$ is a normal subgroup of
$Z$, we have that $Z/(Z\cap\Lambda)$ is a compact group. Hence
$x_0\in \cl{\{z^kx_0:k\in\NN\}}$. Therefore $x_0\in Y=\cl{Hx_0}\setminus
Hx_0$, which is a contradiction. This completes the proof.
\qed

\begin{lemma}  \label{lemma:finite:index}
Let $G$ be a Lie group and $\Gamma$ a closed subgroup of $G$ such that
$G/\Gamma$ has a finite $G$-invariant measure. Let
$\Gamma_1\subset\Gamma$ be a subgroup of finite index in
$\Gamma$. Let $\rho:G/\Gamma_1\to G/\Gamma$ be the natural
quotient map, $x_1=e\Gamma_1$, and $x_0=\rho(x_1)$. Let $H$ be a
closed subgroup of $G$. Then the following statements hold:
\begin{enumerate}
\item[{\rm 1.}]
$\cl{Hx_0}\setminus Hx_0 \mbox{ closed in }G/\Gamma \iff
\cl{Hx_1}\setminus Hx_1 \mbox{ closed in } G/\Gamma_1$.
\item[{\rm 2.}]
$Hx_0 \mbox{ closed in } G/\Gamma \iff
Hx_1 \mbox{ closed in } G/\Gamma_1$.
\end{enumerate} 

\end{lemma}

\proof Let $\Gamma^0$ be the connected component of $e$ in
$\Gamma$. Then $\Gamma^0\subset \Gamma_1$. Let $\Omega$ be a relatively
compact open neighbourhood of $e$ in $H$ such that
$\Omega\inv\Omega\cap \Gamma\subset \Gamma^0$. Then for any
$\gamma,\gamma'\in\Gamma$,
\begin{equation}  \label{eq:disjoint}
\cl{\Omega}\gamma x_1\cap\cl{\Omega}\gamma' x_1\neq\emptyset\implies 
\gamma x_1=\gamma'x_1.
\end{equation}

{\bf [1:$\Rightarrow$]} Since $\Omega x_0$ is open in $Hx_0$, we have that
$\rho\inv(\Omega x_0)\cap \cl{Hx_1}$ is open in $\cl{Hx_1}$. Now
$\displaystyle \rho\inv(\Omega x_0)=\bigcup_{\gamma\in\Gamma}
\Omega\gamma x_1
$.
Therefore by Equation~\ref{eq:disjoint}, $\Omega x_1$ is open in
$\cl{Hx_1}$. Hence $Hx_1$ is open in $\cl{Hx_1}$. This
proves~(1:$\Rightarrow$).

\medskip
{\bf [1:$\Leftarrow$]} There exists $\Gamma_2\subset \Gamma_1$ which is a
normal subgroup of finite index in $\Gamma$. Put $x_2=e\Gamma_2$.
By (1:$\Rightarrow$), $\cl{Hx_2}\setminus Hx_2$ is closed in
$G/\Gamma_2$. Now by Lemma~\ref{lemma:central:project} applied to
$\Gamma_2$ in place of $\Lambda$ and $\Gamma$ in place of $Z$,
we have that $\cl{Hx_0}\setminus Hx_0$ is closed in $G/\Gamma$. This
proves (1:$\Leftarrow$).

\medskip
{\bf [2:$\Leftarrow$]} This holds because $\rho$ is a proper map. 

\medskip
{\bf [2:$\Rightarrow$]} Since $\cl{Hx_0}\setminus Hx_0=\emptyset$ is
closed in $G/\Gamma$, by (1:$\Rightarrow$) $\cl{Hx_2}\setminus Hx_2$
is closed.  By Lemma~\ref{lemma:central:project} applied to $\Gamma_2$ in
place of $\Lambda$ and $\Gamma$ in place of $Z$, we get that
$Hx_2$ is closed. Now (2:$\Rightarrow$) follows from (2:$\Leftarrow$).
\qed

\medskip
Next we recall some of the properties of actions of unipotent subgroup
on homogeneous spaces which will be used in the proof of
Theorem~\ref{nimish:thm:closed:finite}. The properties considered here do not
involve the description of invariant measures and orbit closures for
such actions.

\subsubsection*{\it Nondivergence of unipotent trajectories on finite
volume homogeneous spaces and consequences}

\begin{theorem}[Dani] \label{nimish:thm:Dani:return}
Let $G$ be a connected Lie group and $\Gamma$ a lattice in $G$. Let a
compact set $C \subset G/\Gamma$ and an $\epsilon>0$ be given. Then
there exists a compact set $K\subset G/\Gamma$ such that for any
$\Ad_G$-unipotent element $u\in G$ and any $x\in C$ the following
holds:
\begin{equation} \label{eq:return2}
\frac{1}{N}\sum_{n=1}^N \chi_K(u^nx)>1-\epsilon, \qquad \forall
N\in\NN,
\end{equation}
where $\chi_K$ denotes the characteristic function of $K$ on $G/\Gamma$.
\end{theorem}

\proof For one-parameter unipotent subgroups, the analogous result is
essentially proved in Dani~\cite{Dani:rk=1}; see
\cite[Theorem~6.1]{Dani:Margulis:distribution} for details.  For the
discrete flows, we extend the action of a cyclic unipotent subgroup to
the action of a one-parameter unipotent subgroup as in the beginning
of Section~\ref{sec:extend}. Now the analogous result for the
one-parameter unipotent subgroup action implies the result for action
of a cyclic unipotent subgroup.  \qed

\medskip 
Certain observations due to Margulis in \cite{Margulis:varna}, relating to
Theorem~\ref{nimish:thm:property(D)} and Moore's version of Mautner phenomenon 
lead to the following result.

\newcommand{\cittwo}{\cite[Theorem 2.3]{Shah:distribution}}

\begin{theorem}[\mbox{\cittwo}] \label{nimish:thm:smallest}
Let $G$ be a Lie group, $\Gamma$ a lattice in $G$, and $\pi:G\to
G/\Gamma$ be the natural quotient map. Let $U$ be a
subgroup of $G$ generated by $\Ad_G$-unipotent one-parameter
subgroups. Then there exists the smallest closed connected subgroup $L$
of $G$ containing $U$ such that $\pi(L)$ is closed. Further
$\pi(L)$ admits a finite $L$-invariant measure, which is
$U$-ergodic. 
\end{theorem} 

\begin{corollary} \label{nimish:cor:smallest:density}
Let the notation be as in Theorem~\ref{nimish:thm:smallest}. Then 
the following statements hold.
\begin{enumerate}
\item[{\rm 1.}] 
There exists $g\in L$ such that $\cl{\pi(g Ug\inv)}=\pi(L)$.
\item[{\rm 2.}] 
$\Ad_G(L)\subset\Zcl(\Ad_G(L\cap\Gamma))$. 
\item[{\rm 3.}]
$\pi(N^1_G(L))$ is closed.
\end{enumerate}
\end{corollary}

\proof 
By Theorem~\ref{nimish:thm:smallest}, we have that $U$ acts ergodically on
$\pi(L)$ with respect to a finite $L$-invariant measure. Therefore
statement~1 follows from Hedlund's lemma, statement~2 follows from
statement~1 and Proposition~\ref{nimish:prop:closed:discrete}, and statement~3
follows from statement~1 and Proposition~\ref{nimish:prop:nor}. \qed

\medskip
Next we consider certain properties of actions of subgroups generated
by unipotent elements on finite volume homogeneous spaces of
semisimple groups. 

\begin{proposition} \label{nimish:prop:semisimple}
Let $G$ be a Lie group and $\Gamma$ a lattice in $G$. Suppose that
$G^0$ is a semisimple group with trivial center, $G=\Gamma G^0$, and
$Z_G(G^0)\subset \Gamma$. Let $W\subset G$ be a subgroup such that
$\Zcl(\Ad_G(W))=\Zcl(\Ad_G(\langle \cU\rangle))$, where $\cU$ consists
of $\Ad_G$-unipotent elements of $W$. Then there exists a homomorphism
$\rho:W\to G^0$ such that $\Ad_G(w)=\Ad_G(\rho(w))$ and $wx=\rho(w)x$
for all $x\in G/\Gamma$ and $w\in W$.
\end{proposition} 

\proof Since $G^0$ is semisimple, $\Ad_G(u)\in\Ad_G(G^0)$ for all
$\Ad_G$-unipotent elements of $G$. Since $\Ad_G(G^0)$ is a connected
adjoint semisimple group, it is Zariski closed. Therefore
$\Ad_G(W)\subset\Ad_G(G^0)$. Now since the center of $G^0$ is trivial,
there exists a homomorphism $\rho:W\to G^0$ such that
$\Ad_G(w)=\Ad_G(\rho(w))$ for all $w\in W$. 

Let $w\in W$. Put $\delta=w\inv \rho(w)$. Then $\delta\in
Z_G(G^0)\subset\Gamma$. Now for any $g\in G^0$, 
\[
\rho(w)\pi(g)=w\delta\pi(g)=w\pi(\delta g)=w\pi(g\delta)=w\pi(g).
\]
Since $G/\Gamma=\pi(G^0)$, the above equation holds for all $g\in G$.
\qed

\begin{proposition} \label{nimish:prop:rank2}
Let $G$ be a connected semisimple (real algebraic) group of real rank
$\geq 2$ with trivial center and no nontrivial compact factors. Let
$\Gamma$ be an irreducible lattice in $G$ and $x_0=e\Gamma\in
G/\Gamma$. Let $\{U_j\}$ be a collection of unipotent one-parameter
subgroups of $G$. For each $j$, let $F_j$ be a the smallest closed
connected subgroup of $G$ containing $U_j$ such that $F_jx_0$ is
closed. Let $F$ be the smallest algebraic subgroup of $G$ containing
$F_j$ for all $j$. Then $Fx_0$ is closed, $Fx_0$ has a finite
$F$-invariant measure, and the solvable radical of $F$ is
unipotent. 

In particular, $F$ is the smallest closed connected subgroup of $G$
containing $U_j$ for all $j$ such that $Fx_0$ is closed.
\end{proposition}

\proof By the arithmeticity theorem of Margulis, there exists a
semisimple $\QQ$-group $\tilde G$ and a surjective homomorphism
$\phi:\tilde G\to G$ of real algebraic groups such that $\ker \phi$ is
compact, and $\rho(\tilde G(\ZZ))$ and $\Gamma$ are commensurable
(see~\cite{Zimmer:book}). By Lemma~\ref{lemma:finite:index}, without loss of generality\
we may replace $\Gamma$ by $\rho(\tilde G(\ZZ))$ and assume that
$\Gamma=\rho(\tilde G(\ZZ))$. Let $\bar\phi:\tilde G/\tilde G(\ZZ)\to
G/\Gamma$ be the quotient map associated to $\rho$. Then $\bar\phi$ is
a proper map.

For each $j$ there exists a unipotent one-parameter subgroup $\tilde
U_j$ in $\tilde G$ such that $U_j=\phi(\tilde U_j)$. Let $\tilde
x_0=\bar\phi(x_0)$. Let $\tilde F_j$ be the smallest closed connected
subgroup of $\tilde G$ containing $\tilde U_j$ such that $F_j\tilde
x_0$ is closed. Then by Corollary~\ref{nimish:cor:smallest:density}~(1) and by
\cite[Prop.~3.2]{Shah:distribution}, we have that $\tilde F_j$ is a
$\QQ$-subgroup of $\tilde G$. If $J$ is the smallest $\QQ$-subgroup of
$\tilde G$ containing $\tilde U_j$ then $J\tilde x_0$ is closed. Hence
$\tilde F_j$ is the smallest $\QQ$-subgroup of $\tilde G$ containing
$\tilde U_j$.

Let $\tilde F$ be the smallest algebraic subgroup of $\tilde G$
containing $\tilde F_j$ for all $j$. Then $\tilde F$ is the smallest
algebraic $\QQ$-subgroup of $\tilde G$ containing all $\tilde
U_j$. Therefore the radical of $\tilde F$ is unipotent, and $\tilde
F\tilde x_0$ is closed, and has a finite $\tilde F$-invariant measure.

Since $\bar \phi$ is a proper map, $\phi(\tilde F_j)x_0$ and
$\phi\inv(F_j)\tilde x_0$ are closed. Now since $U_j\subset\phi(F_j)$
and $\phi(\tilde U_j) \subset F_j$, by the minimality, we have that
$F_j=\phi(\tilde F_j)$.

Since $\phi(\tilde F)$ and $\phi\inv(F)$ are algebraic groups, we have
that $\phi(\tilde F)=F$. Since $\bar\phi$ is proper,
$Fx_0=\bar\phi(\tilde F\tilde x_0)$ is closed, and has a finite
$F$-invariant measure. Moreover since $\ker\phi$ is a compact
semisimple Lie group, $$\Rad(F)=\phi(\Rad(\tilde F)).$$ Hence $\Rad(F)$
is unipotent. This completes the proof.  \qed

\subsection*{Proof of Theorem~\ref{nimish:thm:closed:finite}}

The proof is given through a series of claims.

\begin{claim} \label{claim:discrete}
We may assume that $\Gamma$ is a discrete subgroup of $G$, and
$$\Ad_G(G)\subset\Zcl(\Ad_G(\Gamma)).$$ 
\end{claim}

\bproof
Since $G/\Gamma$ admits a finite $G$-invariant measure, by
Proposition~\ref{nimish:prop:closed:discrete}, we have that
$\Ad_G(\cU)\subset \Zcl(\Ad_G(\Gamma))$. Let
$$G_1=\Ad_G\inv(\Zcl(\Ad_G(\Gamma))).$$ Then $W\subset G_1$, $G_1\subset
N_G(\Gamma^0)$, and $\Gamma\subset G_1$. Therefore $G_1/\Gamma$ admits
a finite $G_1$-invariant measure, and hence
$\pi(G_1)$ is closed. Hence replacing $G$ by
$G_1$, we may assume\ that $\Gamma^0$ is normal in $G$. Now
without loss of generality\ we can replace $G$ by $G/\Gamma^0$ and $\Gamma$ by
$\Gamma/\Gamma^0$, and assume that $\Gamma$ is a discrete subgroup of
$G$. Moreover, $\Ad_G(G)\subset \Ad_G(\Zcl(\Gamma))$. \qedclaim
\eproof

Now note that if $F$ is a closed subgroup of $G$ containing $W$ such
that $Fx_0$ is closed, then $WF^0x_0$ is closed, and hence replacing $F$
by $WF^0$ we may assume\ that $F=WF^0$. 

Now suppose $F_i$ ($i=1,2$) are closed subgroups of $G$ such that
$F_i=WF_i^0$ and $F_ix_0$ is closed. Then $Z=F_1^0x_0\cap F_2x_0$ is
an open an closed subset of the closed set $F_1x_0\cap F_2x_0$. Again
if we put $F=F_1\cap F_2$ then $F^0z$ is open in $Z$ for all $z\in
Z$. Therefore $F^0z$ is closed in $Z$ for all $z\in Z$. Hence $Fx_0$
is closed. Moreover $\dim(F^0)< \min(\dim F_1^0, \dim F^0_2)$, unless
$F_1\subset F_2$ or $F_2\subset F_1$. This shows the existence of the
minimal $H$ as assumed in the statement of the theorem.

\begin{claim} \label{claim:G=Gamma G0}
We may assume that $H=H^0(H\cap \Gamma)$.
\end{claim}

\bproof
Take any $u\in\cU$. By Theorem~\ref{nimish:thm:Dani:return}, there exists a
compact set $K\subset G/\Gamma$ such that the set
$\{k\in\NN:u^kx_0\in K\}$ is infinite. Therefore there exists
$k_1,\ldots,k_n\in\NN$ such that $K \cap (\langle u \rangle H^0x_0) \subset 
\cup_{i=1}^n u^{k_i}H^0x_0$. Therefore there exists $k\in\NN$ such that
$u^kx_0\in H^0x_0$.

Put $W_0=\{w\in W: wx_0\in H^0x_0\}$. Then $W_0H^0x_0=H^0x_0$. Put 
\begin{equation} \label{eq:U0}
\cU_0=W_0\cap\left(\cup_{u\in\cU} \langle u\rangle\right).
\end{equation}
Then $\Ad_G(W)\subset\Zcl(\Ad_G(\langle \cU_0\rangle))$. 

Suppose there exists a closed subgroup $L$ of $W^0H^0$ containing
$W^0$ such that $Lx_0$ is closed. Then by Equation~\ref{eq:U0}, we
have that $F=WL$ is a subgroup of $G$. Note that for $w_1, w_2\in W$,
if $w_1Lx_0\cap w_2H^0x_0\neq\emptyset$ then $w_1\inv w_2\in
W^0$. Hence $WLx_0\cap wH^0x_0=wLx_0$ for all $w\in W$. Therefore
$Fx_0$ is closed in $Hx_0$, and hence in $G/\Gamma$. Thus $F=H$, and
$L=H^0$. This shows that replacing $W$ by $W^0$, we may assume\ that
$Hx_0=H^0x_0$. This completes the proof of the claim. \qedclaim    
\eproof

In particular, we may assume\ that $G=WG^0$, and $Wx_0\subset G^0x_0$. Thus
we may assume\ that $G=G^0\Gamma$.  

\begin{claim} \label{claim:smallest}
We may assume that is no proper subgroup $L$ of $G$ containing $W$
such that $Lx_0$ has finite $L$-invariant measure.
\end{claim}

\bproof
Let $L$ be such a subgroup. Then the claim follows from replacement of $G$
by $L$. \qedclaim
\eproof

\subsubsection*{\it Projecting to semisimple factors}

Let $R$ be the connected solvable radical of $G$. Put
$R'=\cl{R\Gamma}^0$. By  Auslander's
theorem~\cite[Theorem~8.24]{Raghunathan:book} $R'$ is solvable. By
Zariski density of $\Gamma$, $R'$ is normal in $G$. Therefore $R=R'$.
Hence $R\Gamma$ is closed. Now since $R$ is normalized by $\Gamma$,
$R\cap\Gamma$ is a lattice in $R$
(see~\cite[Theorem~1.13]{Raghunathan:book}). Therefore by Mostow's
theorem~\cite[Theorem~3.1]{Raghunathan:book}, $R/(R\cap\Gamma)$ is
compact.

Let $\bar C$ be the product of all maximal connected compact normal
subgroups of $G/R$, and let $\bar Z$ be the center of $(G^0/R)/\bar
C$. Since $(G^0/R)/\bar C$ is semisimple and $G=WG^0$, we have that
$\bar Z$ is central in $(G/R)/\bar C$. Put $\bar G=((G/R)/\bar C)/\bar
Z$. Then $\bar G^0$ is a semisimple group with trivial center and no
nontrivial compact normal subgroups. Let $\sigma:G\to\bar G$ be the
natural quotient homomorphism.

Put $\bar\Gamma=\sigma(\Gamma)$. Then $\bar\Gamma$ is a lattice in
$\bar G$, and $\sigma\inv(\bar\Gamma)/\Gamma$ is compact. Let
$\Lambda=\{g\in\bar G:g(\bar G^0\cap \bar\Gamma)g\inv\subset \bar
G^0\cap\bar\Gamma\}$. Then $\Lambda/\bar \Gamma$ is compact. Let
$\bar\sigma:G/\Gamma\to\bar G/\Lambda$ be the natural quotient
map. Then $\bar\sigma$ is a proper map. Put $\bar
x_0=\bar\sigma(x_0)$. 

\begin{claim} \label{claim:project}
It is enough to show that $\sigma(H)^0\bar x_0$ admits a finite
$\sigma(H)$ invariant measure.
\end{claim}

\bproof
Let $\mu_H$ denote a locally finite $H^0$-invariant Borel  measure on
$H^0x_0$. Since $\bar\sigma$ is a proper map, the projected measure
$\bar\mu_H=\bar\sigma_\ast(\mu_H)$ on $\sigma(H^0)\bar x_0$ is a locally
finite and $\sigma(H^0)$-invariant. Now the claim follows from the
uniqueness, up to constant multiple, of the locally finite
$\sigma(H^0)$-invariant measures on $\sigma(H^0)\bar x_0$. \qedclaim
\eproof

Since $Z_{\bar G}(\Lambda\cap \bar G^0)\subset Z_{\bar G}(\Gamma \cap
\bar G^0)\subset\Lambda$, by Proposition~\ref{nimish:prop:semisimple}, there
exists a homomorphism $\rho:W\to\bar G^0$ such that $\Ad_{\bar
G}(\sigma(w))=\Ad_{\bar G}(\rho(w))$ and $\sigma(w)x=\rho(w)x$ for all
$w\in W$ and $x\in \bar G/\bar\Lambda$.

There exist closed normal subgroups $G_1,\ldots,G_k$ of $\bar G^0$
such that $$\bar G^0=G_1\times\dots\times G_k,$$ and if $p_i:\bar G^0\to
G_i$ is the projection on the $i$-th factor, then $\Lambda_i=
p_i(\Lambda)$ is an irreducible lattice in $G_i$. Let $\bar p_i:\bar
G^0/(\bar G^0\cap \Lambda)\to G_i/\Lambda_i$ be the natural quotient
map. Put $\sigma_i= p_i\circ\sigma|_{G^0}$ and $x_i=\bar
p_i(\bar\sigma(x_0))$.  

Fix any $i\in\{1,\dots,k\}$. Let $U$ be the subgroup of $G_i$
generated by the one-parameter unipotent subgroups $\{u(t)\}$
associated to all $u\in\cU$ such that $p_i(\rho(u))=u(1)$. Let $L_i$
be the smallest closed connected subgroup of $G_i$ containing $U$ such
that $L_ix_i$ is closed. Then by Theorem~\ref{nimish:thm:smallest}, $L_ix_i$
has finite $L_i$-invariant measure.

\begin{claim} \label{claim:Li=Gi}
$L_i=G_i$.
\end{claim}

\bproof
Since the fibers of $\bar p_i$ have finite invariant measures, we have
that $p_i\inv(L_i)\bar x_0=\bar p_i(L_ix_i)$ has a finite
$p_i\inv(L_i)$-invariant measure. Therefore $\sigma(W)p_i\inv(L_i)\bar
x_0=p_i\inv(L_i)\bar x_0$ has a finite $\sigma(W)p_i\inv(L_i)$-invariant
measure. Again
$W(p_i\circ\sigma)\inv(L_i)x_0=\sigma\inv(\sigma(W)p_i\inv(L_i))x_0$
has a finite $W(p_i\circ\sigma)\inv(L_i)$-invariant measure. Now the
claim follows from Claim~\ref{claim:smallest}. \qedclaim
\eproof

\begin{claim} \label{claim:radical:compact}
The radical of $H_i=\sigma_i(H^0)$ is compact for all $1\leq i\leq k$.
\end{claim}

\bproof
For any $u\in\cU$, there exists a closed connected subgroup $F_u$ of
$G$ such that $F_ux_0$ has a finite $F_u$-invariant measure, $u\in
N_G(F_u)$, and $\cl{\langle u\rangle x_0}=\langle u\rangle F_u
x_0$. Therefore $F_u\subset H$. 

Now $\sigma_i(F_u)x_i=\bar p_i(\bar\sigma(F_ux_0))$. Hence
$\sigma_i(F_u)x_i$ has a finite $\sigma_i(F_u)$\-invariant
measure. Therefore $\sigma_i(F_u)x_i$ is closed. Let $\{u(t)\}$ be the
one-parameter unipotent subgroup of $G_i$ such that
$p_i(\rho(u))=u(1)$. Since $\Zcl(\langle u(1)\rangle)=\{u(t)\}$, we
have that $\{u(t)\}\subset N_G(\sigma_i(F_u))$ and $$\{u(t)\}\subset
N_{G_i}(H_i).$$  Now
$\cl{\{u(t)\}x_i}=\{u(t)\}\sigma_i(F_u)x_i$. This shows that
\begin{equation}  \label{eq:normalize:sigma(H)}
\{u(t)\}\sigma_i(F_u)\subset N_{G_i}(H_i). 
\end{equation}
\medskip\noindent {\it Case of\/} $\RR$-rank$(G_i)\geq 2$ {\it :} In
this case by Proposition~\ref{nimish:prop:rank2}, the group $L_i$ is the
Zariski closure of the subgroup generated by $\{u(t)\}\sigma_i(F_u)$
for all $u\in\cU$. Therefore by Equation~\ref{eq:normalize:sigma(H)}
and Claim~\ref{claim:Li=Gi}, $H_i$ is a normal subgroup of
$G_i$. Hence the claim holds in this case.

\medskip\noindent{\it Case of\/} $\RR$-rank$(G_i)=1$ {\it :}
Put $F=N_{G_i}(H_i)$. Then $U\subset F$. Now suppose that
the unipotent radical of $F$ is nontrivial. Since $\rank(G_i)=1$, $F$
is contained in a unique minimal parabolic subgroup, say $P$ of
$G$. Let $N$ denote the unipotent radical of $P$. Then $U\subset N$
and $\sigma_i(H)\supset P$. Hence $\{u(t)\}\sigma_i(F_u)\subset P$ for
all $u\in\cU$. Let $M=Z_{G_i}(A)^0$. Since
$$\cl{\{u(t)\}x_i}=\{u(t)\}\sigma_i(F_u)x_i,$$ we have that
$\{u(t)\}\sigma_i(F_u)\subset MN$. Therefore
$MN\cap\Lambda_i\neq\{e\}$. Hence $Nx_i$ is compact. Therefore
$\cl{Ux_i}\subset Nx_i$, and hence $$N=L_i=G_i,$$ which is a
contradiction. Thus $F$ is a reductive subgroup of $G_i$. Since
$\{e\}\neq U\subset F$, we obtain that the solvable radical of
$H_i$ is compact. This completes the proof of the
claim. \qedclaim
\eproof

\bigskip

\noindent {\bf Completion of the proof of the theorem} 
Since $\bar\sigma$ is a proper map, $\sigma(H^0)\bar x_0$ is
closed. By the Claim~\ref{claim:radical:compact}, the radical of
$\sigma(H^0)$ is compact. Now by \cite[Theorem~15]{Margulis:ICM},
$\sigma(H^0)\bar x_0$ has a finite $\sigma(H^0)$-invariant measure. In
view of Claim~\ref{claim:project}, this completes the proof of the
theorem. \qed

\begin{remark} \rm
We may note that the above results in this section are independent of the
classification of ergodic invariant measures and orbit closures for
unipotent flows.
\end{remark}

\subsection*{Limits of ergodic invariant measures for a discrete
unipotent flow}

Our proof of Theorem~\ref{nimish:thm:closure} is based on the following
result. 

\begin{theorem}[Mozes and Shah]  \label{nimish:thm:limit:measure}
Let $G$ be a Lie group and $\Gamma$ a  discrete subgroup of $G$ such
that $G=G^0\Gamma$. Let $u\in G$ be an $\Ad_G$-unipotent element. Let
$\mu_i$ be a sequence of $u$-invariant $u$-ergodic probability
measures on $G/\Gamma$ and such that $\mu_i$ converges to a
probability measure $\mu$ on $G/\Gamma$ in the weak-$\ast$
topology. Suppose that $x_0=e\Gamma/\Gamma\in\supp(\mu)$. Then there
exists a closed subgroup $H$ of $G$ such that the following holds:
\begin{enumerate}
\item[{\rm 1.}] $\mu$ is $L$-invariant and $\supp(\mu)=Lx_0$.
\item[{\rm 2.}] For any sequence $g_i\to e$ in $G$ such that $\cl{\langle
u\rangle g_ix_0}=\supp(\mu_i)$ (such sequences exist due to Hedlund's
lemma), we have
\[
g_i\inv ug_i\in L, \qquad \forall i\GG 0.
\]
(Here $\forall i\GG 0$ stands for the expression `for all $i\geq i_0$, for
some $i_0>0$'.) 

In other words, $u\in L$ and $g_i\supp(\mu)\subset\supp(\mu)$,
$\forall i\GG 0$. 
\item[{\rm 3.}]
$\Ad_G(L)\subset \Zcl(\Ad_G(L\cap\Gamma))$.
\item[{\rm 4.}]
$L=\langle u\rangle L^0$.
\end{enumerate}
\end{theorem}

\proof 
Using the method of Section~\ref{sec:extend}, we extend the action of
$\langle u\rangle$ to the action of a one-parameter unipotent
subgroup. For the action of a one-parameter unipotent subgroup, the
analogous result holds ~\cite[Theorem~1.1]{Mozes:Shah:limit}. From that 
we deduce statements~1 and~2 of the theorem using the intersection with
$\ZZ\cdot G_0$ as in Section~\ref{sec:extend}.

By Proposition~\ref{nimish:prop:closed:discrete}, $\Ad_G(g_i\inv ug_i)\subset
\Zcl(\Ad_G(L\cap\Gamma))$ for all $i\GG 0$. Put $L_1=L\cap
\Ad_G\inv(\Zcl(\Ad_G(L\cap\Gamma)))$. Then $g_i\inv u g_i\in L_1$ for
all $i\GG 0$, and $L_1x_0$ is closed. Therefore
$g_i\inv\supp(\mu_i)=\cl{g_i\inv\langle u\rangle g_i x_0}\subset
L_1x_0$. Hence $\supp(\mu)\subset L_1x_0$. This shows that
$L\subset L_1$, which proves statement~3.  

To obtain the last statement, note that $u\inv g_i\inv u g_i\to
e$. Therefore $u\inv g_i\inv u g_i\in L^0$ for all $i\GG 0$. Put
$L_2=\langle u\rangle L^0$. Since
$L^0x_0$ is closed and $u\in L$, we have that $L_2x_0$ is closed, and
$g_i\inv u g_i \in L_2$ for all $i\GG 0$. Hence $\supp(\mu)\subset
L_2x_0$. This shows that $L_2=L$. This proves statement~4. \qed

\begin{corollary} \label{nimish:cor:limit:measure}
Let $G$, $\Gamma$, $x_0$, and $u$ be as in be as in
Theorem~\ref{nimish:thm:limit:measure}. Suppose further that $\Gamma$ is a
lattice in $G$. Let $g_i\to e$ be any sequence in $G$. Then after
passing to a subsequence, there exists a closed subgroup $L$ of $G$
such that the following holds:
\begin{enumerate}
\item[{\rm 1.}]
$g_i\inv ug_i\subset L$ for all $i>0$. 
\item[{\rm 2.}] 
$Lx_0\subset \cl{\left\{\cup_{i\geq i_0} \langle u\rangle g_i x_0
\right\}}$ for any $i_0>0$.
\item[{\rm 3.}] $\Ad_G(L)\subset \Zcl(\Ad_G(L\cap\Gamma))$ and $L \cap \Gamma$ is a lattice 
in $L$. 
\item[{\rm 4.}] $L=\langle u \rangle L^0$.
\end{enumerate}
\end{corollary}

\proof 
By theorem~\ref{nimish:thm:ratner:closure}, for each $i>0$, the trajectory
$\{u^ng_ix_0:n>0\}$ is uniformly distributed with respect to a
probability measure $\mu_i$ such that $\cl{\langle u\rangle g_i
x_0}=\supp(\mu_i)$. By theorem~\ref{nimish:thm:Dani:return}, given any
$\epsilon>0$ there exists a compact set $K\subset G/\Gamma$ such that
$\mu_i(K)>1-\epsilon$ for all $i>0$. Therefore after passing to a
subsequence, there exists a probability measure $\mu$ on $G/\Gamma$
such that $\mu_i\to \mu$. Now the conclusion of the corollary follows
immediately from Theorem~\ref{nimish:thm:limit:measure}.  \qed

\section{Proof of Theorem~\ref{nimish:thm:closure}}

Let $\pi:G\to G/\Gamma$ be the natural quotient map. Put
$Y=\cl{Wx}$. Let $g\in\pi\inv(x)$. Then replacing $Y$ by $g\inv Y$ and
$W$ by $g\inv W g$, without loss of generality\ we may assume\ that $Y=\cl{\pi(W)}$.

We argue as in the proof of Claim~\ref{claim:discrete} and assume that
$\Gamma$ is discrete.

By Claim \ref{claim:zero}
we may assume that $\cU$ is finite. Using the arguments of the proof
of Claim~\ref{claim:G=Gamma G0}, we may assume\ that $Y\subset
\pi(G^0)$ and $G=WG^0=\Gamma G^0$. 

Now we divide Theorem~\ref{nimish:thm:closure} into the following two
complementary  theorems. 

\begin{theorem} \label{nimish:thm:open:in:closure}
Let the notation and conditions be as above. Then there  exists a
closed subgroup $M$ of $G$ containing $W$ such that $\cl{\pi(W)}$ is
$M$-invariant and $\pi(M)$ is open in $\cl{\pi(W)}$.
\end{theorem}

\begin{theorem} \label{nimish:thm:boundary}
Let the notation be as above. Suppose there exists a closed subgroup 
$M$ of $G$ containing $W$ such that $\cl{\pi(W)}$ is $M$-invariant,
and $\pi(M)$ is open in $\cl{\pi(W)}$. Then $\cl{\pi(W)}=\pi(M)$. Also $W$
acts ergodically with respect to a locally finite $M$-invariant
measure on $\pi(M)$.
\end{theorem}

We will prove Theorem~\ref{nimish:thm:open:in:closure} by closely following
the arguments of the proof of
Theorem~\ref{nimish:thm:measure:locally:finite}. While the proof of
Theorem~\ref{nimish:thm:boundary} involves Theorem~\ref{nimish:thm:limit:measure} on
limits of ergodic invariant measures for unipotent flows as a main new
ingredient. 

\subsection*{Proof of Theorem~\ref{nimish:thm:open:in:closure}}

The proof is given via a series of claims, which are reductions to
special cases, and propositions.

\begin{claim} \label{claims:minimal:closure}
We may assume that $G$ contains no proper closed subgroup $F$
containing $W$ such that $\pi(F)$ is closed.
\end{claim}

\proof By Theorem~\ref{nimish:thm:closed:finite} there exists a smallest closed 
subgroup $F$ of $G$ containing $W$ such that $\pi(F)$ is closed. Moreover, 
$\pi(F^0)$ has a finite $F^0$-invariant measure. As in Claim~\ref{claim:G=Gamma G0}
we may assume that $F=WF^0=F^0(F\cap \Gamma)$. Therefore, without loss of 
generality, we may replace $G$ by $F$. Now the claim follows. \qed


\begin{proposition}  \label{nimish:prop:cap:N}
There exist $H_u\in\cH_u$ for each $u\in\cU$ such that
$$\pi(\capu)$$ contains an open subset of $Y$, say $\Psi$, and  
$\cup_{u\in\cU} \pi(S(H_u,u))$ does not contain any open subset of
$\Psi$.
\end{proposition}

\proof Let $\Omega$ be any nonempty open subset of $\pi\inv(Y)\cap
G^0$. Express $\cU=\{u_1,\ldots,u_k\}$ for some $k\in\NN$. For any
$H\in\cH_{u_1}$, there exist compact sets $\{C_i(H)\}_{i\in\NN}$ such
that $C_i(H)\subset C_{i+1}(H)$ for all $i\in\NN$, and
$$N(H,u_1)=\cup_{i\in\NN} C_i(H).$$ By Theorem~\ref{nimish:thm:ratner:closure},
for every $g\in\Omega$, there exists $H\in\cH_{u_1}$ such that  $g\in
N(H,u_1)$. Hence
\[
\pi(\Omega)\subset\bigcup_{H\in\cH_{u_1}} \bigcup_{i\in\NN}
\pi(C_i(H)).
\]
By Proposition~\ref{nimish:prop:countable}, $\cH_{u_1}$ is countable. Since
$\pi(\Omega)$ is an open subset of $Y$, by Baire's category theorem,
$\pi(C_i(H))$ contains a nonempty open subset of $\pi(\Omega)$ for
some $H \in \cH_{u_1}$ and some $i\in\NN$. Furthermore we can choose
the $H$ with an additional property that for any $F<H$, the set
$\pi(C_j(F))$ does not contain any nonempty open subset of
$\pi(\Omega)$ for any $j\in\NN$. Then there exists a nonempty open
subset $\Omega'$ of $\Omega$ such that
$\pi(\Omega')\subset\pi(N(H,u_1))$ and, due to Baire's category
theorem, $\pi(S(H,u_1))$ does not contain any nonempty open subset of
$\pi(\Omega)$.

Now by Baire's category theorem, there exists $\gamma_1\in\Gamma$ such
that $N(H,u_1)$ contains a nonempty open subset $\Omega_1$ of
$\Omega'$. Put $H_{u_1}=\gamma_1\inv H \gamma_1$. Then
$\Omega_1\subset N(H_{u_1},u_1)$. Since
$\pi(S(H_{u_1},u_1))=\pi(S(H,u_1))$, we have that
$$\pi(S(H_{u_1},u_1))$$ does not contain a nonempty open subset of
$\pi(\Omega_1)$.

For each $2\leq i\leq k$, repeating this procedure for $u_i$ in place
of $u_1$ and $\Omega_1$ in place of $\Omega$, we obtain
$H_i\in\cH_{u_i}$ such that $N(H_{u_i},u_i)$ contains a nonempty open
subset $\Omega_i$ of $\Omega_{i-1}$ and $\pi(S(H_{u_i},u_i))$ does not
contain a nonempty open subset of $\pi(\Omega_i)$.

Thus $\cap_{1\leq i \leq k} N(H_{u_i},u_i)$ contains a nonempty open
subset $\Omega_0=\Omega_k$ of $\Omega$ and, due to Baire's category
theorem, $\cup_{i=1}^k \pi(S(H_{u_i},u_i))$ does not contain a
nonempty open subset of $\Psi=\pi(\Omega_0)$.  \qed

\begin{claim} \label{claims:e:generic}
We may assume that $e\in\capsu$, $\cl{\pi(W)}=Y$, and $\capu$
contains a neighbourhood of $e$ in $\pi\inv(Y)$.
\end{claim}

\proof By Proposition~\ref{nimish:prop:cap:N} there exists $g\in\capsu$ such
that the following holds: $\cl{W\pi(g)}=Y$ and $\capu$ contains a
neighbourhood of $g$ in $\pi\inv(Y)$. If we replace $\Gamma$ by $g
\Gamma g\inv$ and $H_u$ by $g H_u g\inv$, the claim follows.  \qed

\begin{proposition} \label{nimish:prop:L:in:capu}
Let $L$ be a closed connected subgroup of $G$. If $\pi(L)\subset Y$,
then the following holds: $L\subset \capu$ and $L\cap\capsu$ is dense
in $L$.
\end{proposition}

\proof By Claim~\ref{claims:e:generic}, we have that $\capu$ contains
a neighbourhood of $e$ in $L$. Therefore by
Lemma~\ref{lemma:analytic}, we have that $L\subset\capu$.

Now suppose that $\cup_{u\in\cU} S(H_u,u)$ contains a nonempty open
subset of $L$. Then by Baire's category theorem, there exists $u\in\cU$
and $F<H_u$ such that $N(F,u)$ contains a nonempty open subset of
$L$. Therefore by Lemma~\ref{lemma:analytic}, $L\subset N(F,u)$. In
particular $e\in S(H_u,u)$, which contradicts
Claim~\ref{claims:e:generic}. Therefore $\capsu$ contains a dense
subset of $L$.
\qed

\begin{proposition} \label{nimish:prop:Hv:in:capu}
$H_v^0\subset\capu$. 
\end{proposition}

\proof 
By Claim~\ref{claims:e:generic} and
Proposition~\ref{nimish:prop:g:in:N*(H,u)}, for any $v\in\cU$,
$$\pi(H_v^0)\subset\cl{\pi(\langle v \rangle)}\subset Y.$$ Therefore the
present proposition follows from Proposition~\ref{nimish:prop:L:in:capu}. 
\qed

\bigskip\noindent{\bf Definition}
Let $U_u$ and $F_u$ be as defined before the statement of 
Proposition~\ref{nimish:prop:normalize}. Using
Proposition~\ref{nimish:prop:Hv:in:capu}, in place of
Proposition~\ref{nimish:prop:Hv}, we conclude that
Proposition~\ref{nimish:prop:normalize} is valid. Also
Propositions~\ref{nimish:prop:F} and \ref{nimish:prop:W:in:N(F)} are
valid; that is, $F$ is a closed subgroup generated by all $F_u$
($u\in\cU$), $\pi(F)$ is closed and has a finite invariant measure,
and
\[
W\subset N_G^1(F).
\]

\begin{proposition}  \label{nimish:prop:FY=Y}
$Y$ is $F$-invariant.
\end{proposition}

\proof Express $\cU=\{u_1,\ldots,u_k\}$. For $1\leq i \leq k$, put
$L_i=F_{u_1}\cdots F_{u_i}$. Since $e\in N^\ast(H_{u_1},u_1)$, by
Proposition~\ref{nimish:prop:g:in:N*(H,u)}, we have $$\pi(L_1)\subset
\pi(H_{u_1})=\cl{\pi(\langle u_1 \rangle)}\subset Y.$$ Since $F=L_k$,
to prove that $\pi(L_k)\subset Y$ by induction, we assume that
$\pi(L_i)\subset Y$ for some $1\leq i\leq k-1$. Put
$$L_i^\ast=L_i\cap\left(\capsu\right).$$ Then for $f\in L_i^\ast$, by
Proposition~\ref{nimish:prop:g:in:N*(H,u)},
\[
f\pi(F_{u_{i+1}})\subset f\pi(H_{u_{i+1}})=\cl{\langle u_{i+1}
\rangle\pi(f)}\subset Y.
\]
By Proposition~\ref{nimish:prop:L:in:capu}, $L_i^\ast$ is dense in
$L_i$. Hence $\pi(L_{i+1})=\cl{L_i^\ast\pi(F_{u_{i+1}})}\subset
Y$. Therefore $\pi(F)\subset Y$. 

Since $W\in N_G(F)$, we have that 
\[
FY=F\cl{\pi(W)}\subset \cl{\pi(FW)}=\cl{W\pi(F)}\subset \cl{WY}=Y.
\]
\qed

By Proposition~\ref{nimish:prop:nor}, $\pi(N_G^1(F))$ is closed. Also
$W\subset N_G^1(F)$. Therefore by Claim~\ref{claims:minimal:closure} without
loss of generality we may assume that $F$ is a normal subgroup of $G$.

\begin{claim} \label{claims:Z_G(W)}
We may assume that $Z_G(W)$ contains a neighbourhood of
$e$ in $\pi\inv(Y)$. 
\end{claim}

\bproof Due to Proposition~\ref{nimish:prop:FY=Y}, without loss of generality\ we can pass to
the quotient $G/F \Gamma$ and assume that $F=\{e\}$. Now arguing as in the
proof of Claim~\ref{claims:mu(Z_G(W))}, we conclude that $Z_G(W)$
contains a neighbourhood of $e$ in $\capu$. Now since $\capu$ contains
a neighbourhood of $e$ in $\pi\inv(Y)$, the claim follows.
\qed
\eproof

\noindent {\bf Completion of the proof of the theorem}
Define $M_1=\{z\in Z_G(W):\pi(z)\in\pi(W)\}$. Then by
Claim~\ref{claims:Z_G(W)}, we have that the closure $\cl{M_1}$ contains
a neighbourhood of $e$ in $\pi\inv(Y)$. Now for any $z\in M_1$, 
\begin{equation}  \label{eq:Y=zY}
Y=\cl{W\pi(z)}=\cl{z\pi(W)}=zY.
\end{equation}
Therefore if $M$ is the closure of the subgroup generated by $M_1$ and
$W$. Then $MY=Y$ and $\pi(M)$ is open in $Y$. This completes the proof
of the theorem. \qed

\subsection*{Proof of Theorem~\ref{nimish:thm:boundary}}

We intend to prove this theorem by induction on $\dim(G)$. For
$\dim(G)=0$ the theorem is trivial. 

Put $Y=\cl{\pi(W)}$. Without loss of generality we may assume that 
\begin{equation}  \label{eq:M=}
M=\{g\in G: gY=Y\}.
\end{equation}
Put $Y_1=Y\setminus \pi(M)$. Note that $Y_1$ is
$M$-invariant.  We want to show that $Y_1=\emptyset$. Suppose that 
\begin{equation} \label{eq:Y1}
Y_1\neq\emptyset.
\end{equation}

Then arguing as in Proposition~\ref{nimish:prop:cap:N}, for each $u\in\cU$ there
exists $H_u\in\cH$ such that $\cap_{u\in\cU}N(H_u,u)$ contains a
nonempty open subset of $\pi\inv(Y_1)$, say $\Psi$, such that
$\pi(\cup_{u\in\cU}S(H_u,u))$ does not contain any open subset of
$\pi(\Psi)$. Let $g\in \Psi\setminus\cup_{u\in\cU}S(H_u,u)$.
Replacing $\Gamma$ by $g\Gamma g\inv$ and putting $x=\pi(g\inv)$,
without loss of generality\ we may assume\ that $Y=\cl{Wx}$, $Y_1=Y\setminus Mx$, and
\[
e\in \capsu.
\]
Thus $\cl{\pi(\langle u \rangle)}=\pi(H_u)$. 

\emph{Just as in
Claim~\ref{claims:minimal:closure}, we may assume\ that there is no proper
closed subgroup $F$ of $G$ containing $W$ such that $Fx$ is closed.}

In view of (\ref{eq:Y=zY}), (\ref{eq:M=}) and (\ref{eq:Y1}), there
exists a sequence $g_i\to e$ in $G^0\setminus Z_G(W)$ such that
$\pi(g_i)\in Wx$ for all $i\in\NN$. By
Corollary~\ref{nimish:cor:limit:measure}, after passing to a subsequence, we
have the following: For any $u\in\cU$, there exists a closed subgroup
$L_u$ of $G$ containing $u$ such that $\pi(L_u)\subset Y$,
\begin{equation} \label{eq:[g,u]}
g_i\inv u g_i u^{-1} \in L_u^0, \qquad \forall i\GG 0,
\end{equation}
$L_u \cap \Gamma$ is a lattice in $L$, and 
\[
\Ad_G(L_u)\subset\Zcl(\Ad_G(L_u\cap \Gamma)).
\]
Since $L_u=\langle u\rangle L_u^0$ and $\cl{\pi(\langle u
\rangle)}=H_u$, we deduce that
$$L_u\cap\Gamma=(L_u^0\cap\Gamma)(H_u\cap\Gamma).$$ Hence
\begin{equation} \label{eq:Zcl(H^0)}
\Ad_G(L_u)\subset\Zcl(\Ad_G(L_u^0\cap \Gamma )(H_u\cap \Gamma)).
\end{equation}

We claim that $\pi(L_u)\subset Y_1$. Because, if $Mx$ contains
a nonempty open subset of $\pi(L_u)\subset Y$, then $L_u^0\subset M$,
and hence $L_u\subset M$. But $$\pi(M)\cap Mx=\emptyset,$$ which proves
the claim.

Now by arguments as those involved in the proof of
Proposition~\ref{nimish:prop:Hv:in:capu}, 
\begin{equation} \label{eq:H^0}
L_u^0\subset \cap_{v\in\cU}N(H_v,v).
\end{equation}  

For any $v\in \cU$, let $U_v$ denote the subgroup generated by all
one-\-parameter $\Ad_G$-unipotent subgroups of $H_v$, and $F_v$ denote
the closed  connected subgroup of $H_v$ such that
$\cl{\pi(U_v)}=\pi(F_v)$. Then $F=\prod_{v\in\cU} F_v$ is a closed
subgroup of $G$. Also $Fx$ admits a finite $F$-invariant measure, and
$\pi(N_G(F))$ is closed.

\begin{claim} \label{claims:J:normal}
We may assume that $F$ is a normal subgroup of $G$.
\end{claim}

\proof Since $e\in \cap_{v\in\cU} N^\ast(H_v,v)$, by
(\ref{eq:H^0}) and by the arguments as in the proof of
Proposition~\ref{nimish:prop:normalize} we have that
\begin{equation} \label{eq:H0:Lv}
L_u^0\cap \Gamma \subset N_G(H_v^0)\cap N_G(F_v), \forall v\in\cU.
\end{equation}

On the other hand, by the arguments as those involved in
Proposition~\ref{nimish:prop:W:in:N(F)}, we have $v\in N_G(F)$ for all
$v\in \cU$. Since $\pi(N_G(F))$ is closed, we have $H_v \subset N_G(F)$ for 
all $v\in \cU$. Therefore by (\ref{eq:Zcl(H^0)}),
\begin{equation} \label{eq:H:Lu}
L_u\subset N_G(H_u^0) \cap N_G(F).
\end{equation}

Since $g_i\inv u g_i\in L_u$ for all $i\GG 0$, we have that $u\in
g_iN_G(F)g_i\inv$ for all $i\GG 0$. Therefore $W\subset
g_iN_G(F)g_i\inv$ for all $i\GG 0$. For each $i>0$ there  exists
$w_i\in W$ such that $w_ix=\pi(g_i)$. Therefore $W\subset (w_i\inv g_i
N_G(F) g_i\inv w_i)$, and $(w_i\inv g_i N_G(F) g_i\inv w_i)x=w_i\inv
g_i\pi(N_G(F))$ is closed. Therefore by an assumption made earlier in the proof,
$w_i\inv g_i N(F) g_i\inv w_i=G$; or in other words, $F$ is 
a normal subgroup of $G$. \qed

\medskip
Let $\bar G=G/F$ and $\rho:G\to \bar G$ be the quotient
homomorphism. 
Note that $\pi(F)$ is closed. Therefore $\bar\Gamma=\rho(\Gamma)$ is a
discrete subgroup of $\bar G$. Let $\bar\pi: \bar G \to \bar G/\bar
\Gamma$ be the natural quotient. Note that if $\Gamma$ is a lattice in
$G$, then $\bar\Gamma$ is also a lattice in $\bar G$.

\begin{claim} \label{claims:J} 
We may assume that $\rho(H_v^0)$ contains no non\-triv\-ial $\Ad_{\bar
G}$\-uni\-potent one\-parameter subgroup for all $v\in\cU$.
\end{claim}

\proof If it does, then by the same arguments as before on $\bar
G/\bar \Gamma$, we can go modulo another subgroup like $F$ containing
all the one-parameter $\Ad_{\bar G}$-unipotent subgroups of all $H_v$.
\qed

\medskip

By (\ref{eq:H^0}) and the arguments as in the proof of
Claim~\ref{claims:mu(Z_G(W))}, we have 
\begin{equation}  \label{eq:Z(rho(W))}
\rho(L_u^0)\subset Z_{\bar G}(\rho(W)), \qquad \forall u\in\cU. 
\end{equation}
Therefore by (\ref{eq:[g,u]}), we have that $\rho(g_i\inv u
g_i)$ and $\rho(u)$ commute with each other. Therefore $\rho(g_i\inv u
g_i u\inv)$ is an $\Ad_{\bar G}$-unipotent element for all $i\GG
0$. Since $g_i\to e$, for each $u\in\cU$ and each $i\GG 0$, there
exists a unique one-parameter $\Ad_{\bar G}$-unipotent subgroup
$\{u_i(t)\}$ such that 
\[
u_i(1)=\rho(g_i\inv u g_i u\inv).
\]
Therefore $\{u_i(t)\}\subset \rho(L_u^0)$ for all $i\GG 0$. Now since
$\bar\pi(\rho(L_u^0))$ has a finite $\rho(L_u^0)$-invariant measure,
by Proposition~\ref{nimish:prop:closed:discrete} and by
(\ref{eq:H0:Lv}), we have
\[
\{u_i(t)\}\subset N_{\bar G}(\rho(H_v^0)), \qquad
\forall u,v\in \cU,\,\forall i\GG 0.
\]
Now by Claims \ref{claims:J}, Proposition~\ref{nimish:prop:closed:discrete} and 
(\ref{eq:Z(rho(W))}), we have 
\begin{equation} \label{eq:[g,u]:Lv}
\{u_i(t)\} \subset Z_{\bar G}(\rho(H_v)), \qquad
\forall u,v\in\cU, \, \forall i\GG 0.
\end{equation}

Put
\begin{equation}  \label{eq:S:Z(W)}
S_1=\left[\bigcap_{v\in\cU} Z_{\bar G}(\rho(H_v))\right]^0 \subset
Z_{\bar G}(\rho(W)).
\end{equation}
By Lemma~\ref{lemma:centralizer}, $\bar{\pi}(S_1)$ is closed. Also
$\{u_i(t)\} \subset S_1$ for all $u\in\cU$ and for all $i\GG 0$. Let $S$ be
the smallest closed connected subgroup of $S_1$ such that $\bar{\pi}(S)$ is
closed, and $\{u_i(t)\}\subset S$ for all $u\in\cU$ and all
$i\GG 0$. 

\begin{proposition} \label{nimish:prop:S:Zariski}
$\Ad_{\bar G}(S)\subset \Zcl(\Ad_{\bar G}(S\cap \bar\Gamma ))$.
\end{proposition}

\proof Since $\bar\pi(\rho(L_u^0))$ has a finite
$\rho(L_u^0)$-invariant measure, by Theorem\-\ref{nimish:thm:smallest}, for
each $i\GG 0$, there exists a smallest closed connected subgroup
$S_{u,i}$ of $L_u^0$ containing $\{u_i(t)\}$ such that
$\bar\pi(S_{u,i})$ is closed. And by
Corollary~\ref{nimish:cor:smallest:density}(2),
\begin{equation} \label{eq:Sui:Zariski}
\Ad_{\bar G}(S_{u,i})\subset \Zcl(\Ad_{\bar G}(S_{u,i}\cap \bar\Gamma)), \qquad
\forall u\in\cU,\,\forall i\GG 0.
\end{equation}
Since $\bar\pi(S)$ is closed,  $\bar\pi(S\cap
\rho(L_u^0))$ is closed. Therefore by minimality, $$S_{u,i}\subset S
\mbox{~for~ all~} u\in\cU \mbox{~and~ all~} i\GG 0.$$
Put 
\[
S'=S \cap \Ad_{\bar G}\inv \Zcl(\Ad_G(S\cap \bar\Gamma)).
\]
Since $\bar\pi(S)$ is closed, we have that $\bar\pi(S')$ is closed. Since
$\{u_i(t)\}\subset S'$ for all $u\in\cU$ and all $i\GG 0$, by
minimality, $S'=S$. \qed 

\begin{claim}  \label{claims:lattice:S}
We may assume that $S$ is central in $\bar G$. 
\end{claim}

\proof By definition of $S$ and
Corollary~\ref{nimish:cor:smallest:density}(3), $\bar\pi(N^1_{\bar G}(S))$ is
closed.  By~(\ref{eq:S:Z(W)}), we have $\rho(g_i\inv u g_i)\in
N^1_{\bar G}(S)$ for all $u\in\cU$ and $i\GG 0$. Therefore by Zariski
density, $\rho(g_i\inv W g_i)\subset N^1_{\bar G}(S)$ for all
$i\GG 0$. Now by the arguments at the end of the proof of
Claim~\ref{claims:J:normal}, without loss of generality we may assume that $S$ 
is a normal subgroup of $\bar G$. Therefore by 
Proposition~\ref{nimish:prop:S:Zariski} and
Lemma~\ref{lemma:centralizer}, we have $Z_{\bar G}(S)\bar\pi(\rho(g\inv))$ is 
closed. Note that $\rho(W) \subset Z_{\bar G} (S)$ and $x=\pi(g^{-1})$. 
Therefore by an earlier assumption, $\bar G = Z_{\bar G} (S)$. \qed

\begin{proposition} \label{nimish:prop:ST}
The orbit $\bar \pi(S)$ is compact, and there exists a closed subgroup $T$
of $G$ containing $\rho(W)$ such that
$\cl{\bar\pi(\rho(W))}=\bar\pi(T)$. In particular, $\bar\pi(ST)$ is
closed.
\end{proposition}

\proof By Claim~\ref{claims:lattice:S}, we have that $S\cong \RR^k$
for some $k$ and $S\cap \bar\Gamma \cong \ZZ^r$ for some $0\leq r\leq
k$. Then for all $u\in\cU$ and all $i\GG 0$, we have that $S_{u,i}$ is
contained in the subgroup of $S$ corresponding to the $\RR^r$ with respect
to the above isomorphisms. Therefore by the definition of $S$, $S\cong \RR^r$;
that is $\bar\pi(S)$ is compact.

By Theorem~\ref{nimish:thm:open:in:closure}, there exists a closed subgroup
$T$ of $\bar G$ containing $\rho(W)$ such that $\cl{\bar\pi(\rho(W))}$
is $T$-invariant and $\bar\pi(T)$ is open in $\cl{\bar\pi(\rho(W))}$.
 
Since $g_i\not\in Z_G(W)$ for all $i>0$, we have that $\dim \bar
G/S<\dim G$. Therefore by the induction hypothesis applied to 
$\bar G/S\bar \Gamma$, and
Lemma~\ref{lemma:central:project}, we conclude that $\bar\pi(T)$ is
closed; that is, $\cl{\bar\pi(\rho(W))}=\bar\pi(T)$. This completes
the proof of the proposition.\qed

\bigskip

\noindent {\bf Completion of the proof of the theorem}
Since $\rho(g_i\inv u g_i u\inv)\in S$ for all $u\in \cU$ and all
$i\GG 0$ and $S$ is normal in $\bar G$ $\rho(g_i^{-1} w g_i w^{-1}) \in S$ for all $w\in\langle \cU \rangle$. Since $\Ad_{\bar G}(\rho(W))\subset
\Zcl(\langle \Ad_{\bar G}(\cU)\rangle)$, and since $S$ is a normal subgroup of 
$\bar G$, we deduce that $\rho(g_i\inv w g_i w^{-1})\in S$ for all $w\in W$. This 
shows that $\rho(g_i\inv W g_i)\subset ST$ for all $i\GG 0$. Hence by the
arguments as in the last part of the proof of Claims~\ref{claims:J:normal},
we conclude that $\bar G=ST$.

Since $g_i\in ST$, $S\subset Z(\bar G)$ and $\rho(W)\subset T$, we
have that $$\rho(g_i\inv w g_i w\inv)\in T \mbox{~for~ all~} w\in W \mbox{
~and~ all~} i\GG 0.$$ Therefore by the definition of $S$, we have $S\subset
T$. Hence $T=\bar G=G/F$. Now $\cl{\bar\pi(\rho(W))}=\pi(T)$ and
$\cl{\pi(W)}$ is invariant under $F=\ker \rho$. Hence
$\cl{\pi(W)}=\pi(G)$. Since $Y_1\supset\cl{\pi(W)}$ and $Y_1\cap
Mx=\emptyset$, we have a contradiction. Hence $Y_1=\emptyset$. Thus
$\cl{Wx}=Mx$. 

Now using the arguments as in the proof of
Theorem~\ref{nimish:thm:measure:locally:finite}, it is straightforward to
verify that $W$ acts ergodically with respect to a locally finite
$M$-invariant measure on $Mx$. This completes the proof of the
theorem. \qed

\section{Some Consequences}

In this section we derive the corollaries, which are stated in the
introduction, of the descriptions of invariant measures and orbit
closures for the actions of subgroups generated by unipotent elements.

\subsection*{Proof of Corollary~\ref{nimish:cor:lattice}}

The proof given below is due to Dave Witte.

The proof for the description of the ergodic invariant measures follows
from Theorem~\ref{nimish:thm:measure} and the arguments as in \cite[Proof of
Corollary~5.8]{Witte:quotients}.

\subsubsection*{Proof for the description of orbit closures} 

To describe the orbit closures without loss of generality\ we may assume\ that $$x=e\Gamma \mbox{~and~}
G=\cl{W\Gamma}=\Gamma G^0.$$ Now as in Claim~\ref{claim:discrete},
without loss of generality\ we may assume\ that $\Gamma$ is discrete. 

Let $G'=W\times G$, $\Gamma'=\Lambda \times \Gamma$, and $\Delta:W\to
G'$ be the diagonal embedding. Note that $G'/\Gamma'=W/\Lambda\times
G/\Gamma$. Let $$x_1=e\Lambda\in W/\Lambda, \ x_2=e\Gamma\in G/\Gamma,
\mbox{~and~} x'=(x_1,x_2).$$ Note that $G'$, $\Gamma'$, and $\Delta(W)$ satisfy
the conditions of Theorem~\ref{nimish:thm:closure}; and hence there exists a
closed subgroup $F'$ of $G'$ containing $\Delta(W)$ such that
$\cl{\Delta(W)x'}=F'x'$. 

We claim that 
\begin{equation} \label{eq:Delta(Lambda)}
(x_1,\cl{\Lambda x_2})=\cl{\Delta(W)x'}\cap (x_1,G/\Gamma).
\end{equation}

Since $\cl{\Delta(\Lambda)x'}=(x_1,\cl{\Lambda x_2})$, we have that
$(x_1,\cl{\Lambda x_2})\subset \cl{\Delta(W)x'}\cap(x_1,G/\Gamma)$. To
prove the opposite inclusion, suppose that $(x_1,x)\in
\cl{\Delta(W)x'}$. Let $\{w_n\}\subset W$ with $w_nx_1\to x_1$ and
$w_nx_2\to x$. Because $w_nx_1\to x_1$, there exist sequences
$\{\lambda_n\}\subset \Lambda$ and $\delta_n\to e$ in $W$ such that
$w_n=\delta_n\lambda_n$ for all $n\in\NN$. Therefore $\lambda_n x_2 =
\delta_n\inv w_n x_2\to x$. Thus $x\in \cl{\Lambda x_2}$. This proves
the claim.

Now by Equation~\ref{eq:Delta(Lambda)}, 
\begin{eqnarray*}
(x_1,\cl{\Lambda x_2})
&=& F'x' \cap (x_1, Gx_2)\\
&=& (F'\cap (\Lambda \times G))x'\\
&=& (x_1, L x_2),
\end{eqnarray*}
where $L$ is the projection of $F'\cap (\Lambda \times G)$ into
$G$. Thus $\cl{\Lambda x_2}=Lx_2$. Since $Lx_2$ is closed, replacing $L$ by 
$\bar L$ we have that $L$ is a closed subgroup of $G$.

By Theorem~\ref{nimish:thm:closure}, $(F')^0x'$ has a finite invariant
measure. Let $F'_0$ be any subgroup of $F'$ containing $(F')^0$ such
that $F'_0x'$ has a finite $F'_0$-invariant measure.  Because the
stabilizer of $x'=(x_1,x_2)$ in $F'_0$ is contained in $F'_0\cap
(\Lambda \times G)$, this implies that $(F'_0\cap (\Lambda\times
G))x'$ also has a finite invariant measure (see \cite[Lemma~1.6,
p.~20]{Raghunathan:book}). Thus, letting $L'$ be the projection of
$F'_0\cap (\Lambda \times G)$ into $G$, we see that $L'x_2$ has a
finite invariant measure. Because $L'$ is open in $L$, we have that
$L^0$ is the identity component of $L'$. Hence we conclude that
$L^0x_2$ has a finite invariant measure. This completes the main part
of the proof of the corollary. 

Note that if $F'x'$ has a finite $F'$-invariant measure, then
$L'=L$. Therefore $Lx_2$ has a finite $L$-invariant measure. Note that
if $W$ is connected, then $F'x'$ is connected, and
hence $Lx_2$ has a finite $L$-invariant measure.
\qed

\subsection*{Proof of Corollary~\ref{nimish:cor:Lambda.g.Gamma}}

First without loss of generality\ we may assume\ that $\Lambda$ is irreducible.  Since $G$ is a
connected semisimple group with no nontrivial compact factors, $G$ is
generated by $\Ad_G$-unipotent one-parameter subgroups. By
Corollary~\ref{nimish:cor:lattice}, applied to $W=G$, there exists a closed
subgroup $F$ of $G$ containing $\Lambda$ such that
$\cl{\Lambda\Gamma}=F\Gamma$. Therefore $F^0$ is normalized by
$\Lambda$. By Borel's density theorem, $\Ad_G(\Lambda)$ is Zariski
dense in $\Ad_G(G)$. Therefore $F^0$ is a normal subgroup of $G$.

Note that $\Lambda F^0$ is an open, and hence a closed, subgroup of
$F$. Therefore the projection of $\Lambda$ on $G/F^0$ is discrete. Now
since $\Lambda$ is an irreducible lattice in $G$, either $F^0=\{e\}$
or $F=G$. Thus either $\Lambda\Gamma$ is discrete, or $\Lambda\Gamma$
is dense in $G$. Since $W=G$ is connected, by
Corollary~\ref{nimish:cor:lattice}, $F\Gamma/\Gamma$ has a finite
$F$-invariant measure. Therefore if $F^0=e$ then
$F\Gamma/\Gamma=\Lambda\Gamma/\Gamma$ is finite. This shows that
$\Lambda\cap\Gamma$ is of finite index in $\Lambda$. Therefore
$\Lambda\cap\Gamma$ is a lattice in $G$, and hence it is a subgroup of
finite index in $\Gamma$.  \qed

\bigskip\noindent{\it Acknowledgment.}  The author would like to thank
Dave Witte for fruitful discussions, and for pointing out several
inaccuracies in an earlier version of this article. The proof of
Corollary~\ref{nimish:cor:lattice} in this article is due to Dave Witte. The
research reported in this article was carried out at the Tata Institute of
Fundamental Research, Bombay, and at the Institute for Advanced Study
(IAS), Princeton. The author acknowledges the support at the IAS
through the NSF grant DMS~9304580.

\vskip1cm

\begin{flushleft}
Institute for Advanced Study, Princeton, NJ 08540, USA.

\vskip 3mm

{\it Permanent Address:}\\
Nimish A. Shah \\
School of Mathematics\\
Tata Institute of Fundamental Research \\
Homi Bhabha Road\\
Mumbai 400005 \\
India\\[3mm]
\emph{E-mail : nimish@math.tifr.res.in}
\end{flushleft}


\begin{thebibliography}{Ab1}

\bibitem[D1]{Dani:invent78} 
S.G. Dani, 
\emph{Invariant measures of horospherical flows on noncompact
homogeneous spaces},
Invent. Math. {\bf 47} (1978), 101--138. 

\bibitem[D2]{Dani:invent2} S.G. Dani,  
\emph{Invariant measures and minimal sets of horospherical flows},
Invent. Math. {\bf 64} (1981), 357--385.

\bibitem[D3]{Dani:version:borel:density}
S.G. Dani,
\emph{On ergodic quasi-invariant measures of group automorphism}, 
Israel J. Math. {\bf 43} (1982), 62--74.

\bibitem[D4]{Dani:rk=1} 
S.G. Dani,  
\emph{On orbits of unipotent flows on homogeneous spaces II}, 
Ergod. Th. and Dynam. Sys. {\bf 4} (1984), 25--34.

\bibitem[DM]{Dani:Margulis:distribution} 
S.G. Dani and G.A. Margulis,  
\emph{Limit distributions of orbits of unipotent
flows and values of quadratic forms},  
Adv. Soviet Math. (publ. by the Amer. Math. Soc.) {\bf 16} (1993), 91--137. 

\bibitem[M1]{Margulis:varna} 
G.A. Margulis,
\emph{Lie groups and ergodic theory}, 
in: Algebra --- Some Current Trends, (Proceedings, Varna 1986), Ed: L.L. Avramov 
and others, pp. 130--146, Lect. Notes Math. {\bf 1352}, Springer-Verlag, 1988.


\bibitem[M2]{Margulis:selberg} 
G.A. Margulis,
\emph{Discrete subgroups and ergodic theory}, 
in: Number Theory, Trace Formulas and Discrete
Subgroups, Symposium in honor of A. Selberg, Oslo 1987, pp. 377--398, Academic
Press, 1989.

\bibitem[M3]{Margulis:ICM}
G.A.~Margulis, 
\emph{Dynamical and ergodic properties of subgroups  actions on
  homogeneous spaces with applications to number theory}, 
in: Proceed. of the Internat. Congr. of Math. (Kyoto, 1990), pp. 193--215,
Mathematical Society of Japan, Tokyo, 1991.

\bibitem[Mo]{Moore:Mautner}
C.C.~Moore,
\emph{The Mautner phenomenon for general unitary representations},
Pacific J.\ Math. {\bf 86} (1980), 115--169.

\bibitem[MS]{Mozes:Shah:limit}
S. Mozes and N. Shah,
\emph{On the space of ergodic invariant measures of unipotent flows},
Ergod. Th. Dynam. and  Sys. {\bf 15} (1995),
149--159.

\bibitem[R]{Raghunathan:book} 
M.S. Raghunathan,  \emph{Discrete Subgroups of Lie Groups}, Springer-Verlag, 1972.

\bibitem[Ra1]{Ratner:solvable}
M.~Ratner, \emph{Strict measure rigidity for unipotent subgroups of
solvable groups}, Invent. Math. {\bf 101} (1990), 449--482.

\bibitem[Ra2]{Ratner:measure}
M.~Ratner, \emph{On Raghunathan's measure conjecture},
Ann. Math. {\bf 134} (1991), 545--607.

\bibitem[Ra3]{Ratner:distribution}
M.~Ratner, \emph{Raghunathan's topological conjecture and distributions of 
unipotent flows}, Duke Math. J. {\bf 63} (1991), 235--290.

\bibitem[Ra4]{Ratner:ICM} 
M.~Ratner,
\emph{Interactions between ergodic theory, Lie groups, and number
theory}, 
in: Proceed. of the Internat. Congr. Math. (Zurich, 1994), pp. 157--182,
Birkh\"auser, 1995. 

\bibitem[Sh]{Shah:distribution}
N.A. Shah,
\emph{Uniformly distributed orbits of certain flows on homogeneous
spaces}, Math. Ann. {\bf 289} (1991), 315--334.

\bibitem[W]{Witte:quotients} 
D. Witte, 
\emph{Measurable quotients of unipotent translations on
homogeneous spaces},  
Trans.\ Amer.\ Math.\ Soc. {\bf 345} (1994),
577--594. 

\bibitem[Z]{Zimmer:book} 
R.J. Zimmer,
\emph{Ergodic Theory and Semisimple Groups}, Birkh\"auser, 1984. 


\end{thebibliography}
\end{document}